\documentclass[preprint,12pt]{article}
\usepackage{amssymb}
\usepackage{mathrsfs}
\usepackage{amsfonts}
\usepackage{graphicx}
\usepackage{enumerate}
\usepackage{indentfirst}
\usepackage{colortbl,dcolumn}
\usepackage{color}
\usepackage{amsmath}
\usepackage{psfrag}
\usepackage{booktabs}
\usepackage{comment}
\usepackage{cite}
\usepackage{amsthm}
\numberwithin{equation}{section}
\usepackage[top=1in, bottom=1in, left=0.8in, right=0.8in]{geometry}
\newcommand{\R}{\mathbb{R}}
\newcommand{\N}{\mathbb{N}}
\newcommand{\E}{\mathbb{E}}
\renewcommand{\P}{\mathbb{P}}

\newcommand{\dd}{\text{d}}
\newtheorem{thm}{Theorem}[section]

\newtheorem{lem}[thm]{Lemma}

\newtheorem{example}[thm]{Example}
\newtheorem{assumption}[thm]{Assumption}
\begin{document}
\title{
	Strong convergence rates for long-time approximations of 
	SDEs with non-globally Lipschitz continuous coefficients
\footnotemark[2] 
\footnotetext[2]{This work was supported by Natural Science Foundation of China (12071488, 12371417, 11971488). The authors want to thank Prof. Zhihui Liu and Yajie She for helpful comments.
}}

\author{Xiaoming Wu$^{\, \text{a,b}}$, Xiaojie Wang$^{\, \text{b}}$ 
\thanks{Corresponding author.
	\newline E-mail addresses:
12331004@mail.sustech.edu.cn,x.j.wang7@csu.edu.cn.} 
\\
\footnotesize $^\text{a}$ Department of Mathematics , Southern University of Science and Technology, Shenzhen, Guangdong, P. R. China
\\
\footnotesize $^\text{b}$ School of Mathematics and Statistics, HNP-LAMA, Central South University, Changsha, Hunan, P. R. China
}
\maketitle	
\begin{abstract}
This paper is concerned with 
long-time strong approximations of SDEs with 
non-globally Lipschitz coefficients.
Under certain non-globally Lipschitz conditions, 
a long-time version of fundamental strong convergence theorem 
is established for general one-step time discretization schemes. 
With the aid of the fundamental strong convergence theorem, 
we prove the expected strong convergence rate over infinite time for two 
types of schemes such as the backward Euler method and 
the projected Euler method in non-globally Lipschitz settings. 
Numerical examples are finally reported to confirm 
our findings.	
	\par 
	{\bf Keywords:} SDEs with non-globally Lipschitz coefficients,
 backward Euler method, projectd Euler method, strong convergence 
\end{abstract}

\section{Introduction}
\noindent
Stochastic differential equations (SDEs) find prominent applications in engineering, physics, chemistry, finance and many other branches of science.
Nevertheless,
nonlinear SDEs rarely have closed-form solution available and one often resorts to 
their numerical approximations.
To numerically study SDEs,
the globally Lipschitz conditions
are imposed on the coefficient functions of SDEs \cite{kloeden1992stochastic,milstein2004stochastic}.
Under the restrictive assumptions, fundamental mean-square convergence theorems were established for one-step approximation schemes of SDEs on the finite time horizon \cite{milstein2004stochastic} and the infinite time horizon \cite{li2022sqrt}.
However,
the coefficients of most nonlinear SDEs from applications violate the traditional but 
restrictive assumptions.
What happens when the commonly used 
schemes in the global Lipschitz settings are applied to SDEs with non-globally Lipschitz coefficients?
As asserted by \cite{hutzenthaler2011strong},
the most popular Euler-Maruyama method produces divergent 
numerical approximations when used to solve SDEs 
with super-linearly growing coefficients over finite interval $[0, T]$.
In the literature, a large amount of attention 
has been attracted to construct and analyze 
convergent numerical approximations of SDEs with super-linearly 
growing coefficients,
see, e.g., 
\cite{beyn2016stochastic,hutzenthaler2015numerical,hutzenthaler2012strong,mao2015truncated, tretyakov2013fundamental,andersson2017mean,beyn2017stochastic,wang2023mean,gan2020tamed,higham2002strong,hutzenthaler2020perturbation,hutzenthaler2018exponential,kelly2023strong,mao2013strong,mao2013strong1,sabanis2016euler,guo2018truncated,kelly2018adaptive,wang2020mean,fang2020adaptive,wang2013tamed,neuenkirch2014first}.
%
%\cite{kloeden1992stochastic,milstein2004stochastic} 
%
Indeed, a majority of existing works are devoted to 
developing and analyzing schemes for strong approximations of SDEs with non-globally Lipschitz conditions over the finite time interval $[0, T]$. 
%
\begin{comment}
  As implied by \cite{giles2008multilevel}, strong approximations are of particular importance for
the computation of statistical quantities of the solution process of SDEs through computationally efficient multilevel Monte Carlo (MLMC) methods.  
\end{comment}
%
In the strong convergence analysis over  the finite time interval $[0, T]$, a powerful tool is the fundamental strong convergence theorem due to \cite{milstein2004stochastic} in the globally Lipschitz setting and its counterpart in a non-globally Lipschitz setting \cite{tretyakov2013fundamental}.
%where the approximation errors are measured in the form of 
%$ \left\| X(T) - Z(T) \right\|  _ { L ^ q ( \Omega; \R ^ { d \times m } ) } $. 
%Here $ X(T) $ and $ Z(T) $ 
%denote the exact solution and the numerical solution of SDEs, 
%respectively.
%In particular, for solving the problem that the additional computational effort for implementation of the backward Euler method,  Hutzenthaler and et.al \cite{hutzenthaler2012strong} first proposed the tamed Euler method and obtained the strong  convergence rate $ \frac{1}{2} $ in non-globally Lipschitz conditions.  Besides, for obtaining better outcomes, some new numerical methods were utlized to study the strong approximations \cite{liu2013strong,mao2015truncated,guo2018truncated,kelly2018adaptive,wang2020mean}. 
%It is worth noting that these numerical methods are based on finite  time to analyze the strong approximations of SDEs with  non-globally Lipschitz assumptions. 

However,
the long-time approximations of SDEs also plays an important role in many scientific areas 
such as high-dimensional sampling, Bayesian inference, statistical physics and machine learning \cite{dalalyan2017theoretical,song2021score,welling2011bayesian,hong2019invariant}.
In the last decades, there are some important works devoted to long-time approximations of SDEs with globally Lipschitz continuous coefficients, see, e.g., \cite{talay1990second,talay2002stochastic,mattingly2010convergence} and references therein.
More recently, some researchers examined long-time approximations of SDEs in non-globally Lipschitz setting
\cite{brosse2019tamed,li2022sqrt,chen2022stochastic,liu2023backward,li2019explicit,fang2020adaptive,liu2023strong,mattingly2002ergodicity,brehier2023approximation,pang2024linear}, to just mention a few.

%Fang and Giles \cite{fang2020adaptive} introduced an adaptive timestep construction for an Euler-Maruyama approximation of ergodic SDEs with non-globally Lipschitz drift over an infinite time interval, because of the adaptive $h_n$  is a function of the current approximate $Z_n$, according to the strategy, one can solve the problem that the uniform timestep explicit Euler-Maruyama discretization causes the scheme not be stable.
%For some small enough $p$, Liu and Wang \cite{liu2023strong} proved the backward Euler method is convergent in the infinite horizon for SDEs with superlinear drift and diffusion coefficients, where the coefficient of the linear term in the drift coefficient is allowed to be positive. 
%As we all known, Monte Carlo methods are a very general and useful approach for the estimation of expectations. However, they can be computationally expensive or inefficient in high dimensions. But for  ergodic SDEs, these problems can be solved by multilevel Monte Carlo.  Thus the strong convergence over infinite time becomes very significant.

To the best of our knowledge, the fundamental strong convergence theorem for long-time approximations of 
	SDEs with non-globally Lipschitz continuous coefficients is still absent in the literature. 
In this paper, we attempt to fill the gap by establishing a long-time version of fundamental strong convergence theorem for general one-step time discretization schemes applied to SDEs with non-globally Lipschitz coefficients (Theorem \ref{thm:global-error}). 
As applications of the fundamental strong convergence theorem, 
we prove the expected strong convergence rate 
of order $\tfrac12$ over infinite time for two 
types of time-stepping schemes such as the backward Euler method and 
the projected Euler method in non-globally Lipschitz settings (Theorem \ref{thm:global-error-backward} and Theorem \ref{thm:global-error-projected}).
\begin{comment}
   Moreover, applications to the Multilevel Monte Carlo method for invariant distribution of SDEs are also examined (see Theorem \ref{thm:MLMC for invariant measure}). 
\end{comment}

%
%Our contribution is two-folded:
%\begin{itemize}
%    \item 
%    Theorem \ref{thm:global-error} can be employed under infinite time compared to Theorem $ 2.1 $ of \cite{tretyakov2013fundamental} ;
%    \item 
%   The strong error analysis for long-time approximations under a relaxed condition compared to \cite{li2022sqrt}.
%\end{itemize}

The outline of the paper is as follows.
In the next section, 
we present some standard notations
and assumptions that will be used in our proofs and obtain the long-time fundamental strong convergence theorem.
In sections \ref{sec:backward} and \ref{sec:projected}, 
we consider the long-time strong convergence analysis of two well-known methods,
e.g., 
the backward Euler scheme and
the projected Euler scheme in the non-globally Lipschitz settings.
%In section \ref{sec:invariant distribution}, we give an application to multilevel Monte Carlo for invariant distribution.
Numerical experiments are  presented in section \ref{sec:numerical experiments}
to verify our theoretical findings.

\section{The long-time fundamental strong convergence theorem}
 \label{sec:settings-stong-convergence-theorem} 
\noindent
\textbf{Notation}: Let $ ( \Omega, \mathcal{ F }, \P ) $ be a complete probability space and
 $ \{ \mathbf{ \mathcal{ F } } _ { t } ^ { w } \}_{0 \leq t  < \infty}$ be an increasing family subalgebras induced by $w ( t )$ for 
$ 0 \leq t  < \infty$, 
where $( w ( t ),\mathbf{ \mathcal{ F } } _ { t } ^ { w })=( ( w _ 1 ( t ), w _ 2 ( t ), \cdots, w _ m ( t ) ) ^ { T }, \mathbf{ \mathcal{ F } } _ { t } ^ { w })$
is an $m$-dimensional Wiener process.
Let $ | \cdot |$ and
 $ \left\langle \cdot , \cdot \right\rangle $
be the Euclidean norm and the inner product of vectors in
 $\mathbb{ R }^ {d}$,
respectively.
By
$ M ^ {T} $ we denote the transpose of a vector or matrix $M$.
For a given  matrix $ M $,
we use
$ \left\| M \right\| := \sqrt{ trace( M ^ T M ) }$
to denote the trace norm of $M$.
On the probability space
    $ ( \Omega, \mathcal{ F }, \P ) $,
we use $ \E $ to denote the expectation and
 $ L ^ q ( \Omega; \R ^ { d \times m } )$,
  $ q \in \N $,
to denote the family of
  $ \R ^ { d \times m } $-valued variables with the norm defined by
  $ \left\| \xi \right\| _ 
     { L ^ q ( \Omega; \R ^ { d \times m } ) }
   =
     ( \E [ \left\| \xi \right\| ^ q ] ) ^ { \frac{1}{q} }
      < \infty $.
      And $\lfloor x \rfloor$ denotes the integer part of $x$.
      
Consider the autonomous SDEs in the It\^{o} form of 
\begin{equation}
 \label{eq:form-SDE}
  \left\{
   \begin{aligned}
	\dd X ( t )
   &
   =
     f ( X ( t ) )
      \, \dd t + 
       g ( X ( t ) )
        \, \dd W ( t ),
   \quad
     t \geq 0,
   \\
	X ( 0 )
   &
   =
    X _ { 0 }
     \in 
      \mathbb{ R } ^ d,
	\end{aligned}
     \right.
\end{equation}
where
 $f = ( f ^ { 1 }, f ^ { 2 }, \cdots, f ^ { d } ) ^ { T }$ 
is the drift coefficient function and
 $g = ( g ^ { i j } ) _ { d \times m } = ( g ^ { 1 }, g ^ { 2 }, 
  \cdots, g ^ { m } ) = ( g _ { 1 } ^ { T }, g _ { 2 } ^ { T }, 
   \cdots, g _ { d } ^ { T } ) ^ { T }$
is the diffusion coefficient function.
It is worth noting that $f ^ { i }$ is a real-valued function,
  $g ^ { i }$ is a $d$-column vector function and $g _ { j }$ is a $m$-row vector function.
In addition,
 $\{ W _ t \} _ { t \geq 0 }$
is an $m$-dimensional Wiener process and the initial data 
${ X _ 0 }$ is independent of $ w$.

We aim to establish a fundamental strong convergence theorem 
for general one-step approximation schemes, used to approximate SDEs \eqref{eq:form-SDE}, 
with uniform step size $h=\tfrac{1}{N}$, $ N \geq 1$.
Here we introduce a new notation $ X ( t, x; s) $ 
for $ t \leq s < \infty $, 
which denotes the solution of \eqref{eq:form-SDE} with the 
initial condition $ X ( t, x; t) = x $. 
When we write $ X ( t ) $, $ t \geq 0$, 
we mean a solution of SDEs \eqref{eq:form-SDE} 
with the initial 
value $ X (0) = X _ 0 $. 
In addition, 
we introduce  one-step approximation
 $ Z ( t, x; t + h ) $ 
 for 
 $ X ( t, x; t+h )$,
where $ t \geq 0, 0<h<1$, 
which is defined as follows:
\begin{equation}
	\label{eq:Z_t-def}
	Z( t, x; t + h ) 
	=
	x
	+
	\Psi ( t, x, h;
	\xi _ t),
\end{equation}
where $ \Psi $ is a function from
$[ 0, \infty ) \times \mathbb{ R } ^ { d }
\times ( 0, T ] \times\mathbb{ R } ^ { m } $ to $ \mathbb{ R } ^ { m }$, 
$ \xi _t $ is a random variable
with sufficiently high-order moments. 
Moreover, 
we define 
 $ Z( t, x; t + h ) $ as an approximation of the solution 
at $ t + h $ with initial value $ Z ( t, x; t) = x $. 
Then, 
a numerical approximation
$ \{ Z _ k \} _ {  k \geq 0 } $
can be constructed
on the uniform mesh grid
$\{ t _ k = kh, k \geq 0 \} $,
given by
\begin{equation}
	\label{eq:Z_k-def}
	Z _ 0 
	=
	X _ 0,
	Z _ { k + 1 }
	= 
	Z _ k
	+
	\Psi
	( t _ k, Z _ k, h; \xi _ k ),
	k \geq 0,
\end{equation}
where $ \xi _ k $ for $ k\geq 0 $
are independent of
$ Z _ 0, Z _ 1, \cdots, Z _ { k - 1 }$,
$\xi _ 0, \xi _ 1, \cdots, \xi _ { k - 1 } $.
Alternatively,
one can write 
\begin{equation}
	\label{eq:Z_k-form}
	Z _ 0
	= 
	X _ 0,
        Z _ { k + 1 }
	= 
	Z ( t _ k, Z _ k; t _ { k + 1 } )
	= 
	Z ( t _ 0, Z _ 0; t _ { k + 1 } ),
	k \geq 0.
\end{equation} 
To facilitate the strong analysis, we impose the following assumptions.
 \begin{assumption}
 	\label{ass:SDEs-xy-contractive}
  Suppose that the drift and diffusion coefficients of SDEs 
  given by \eqref{eq:form-SDE} 
  satisfy a contractive monotone condition.
  Additionally,
there exists a non-integer constant
$ p ^ { * } \in [1 , \infty )$ 
and   
$ \alpha _ 1  \in ( 0, \infty )$  
such that 
\begin{equation}
 \label{eq:SDEs-xy-contractive}
  \left \langle
    x - y,  f ( x ) - f ( y )
     \right\rangle
   +
     \tfrac{ 2p ^ { * } - 1 } { 2 } 
  	   \left\|
  	    g ( x ) - g ( y )
  	     \right\| ^ { 2 } 
   \leq
     -\alpha_1
       |  x - y | ^ { 2 }.
\end{equation}
Assume that f satisfies the polynomial growth Lipschitz condition, 
more accurately,
there exist positive constants 
$ \kappa \geq 1 $ and $ c _ 1 >0 $ 
such that
\begin{equation}
	\label{eq:f-xy-condition}
	|  f ( x ) - f ( y ) | ^ { 2 } 
	\leq
	c _ 1
	\big(
	1 +
	|  x  | ^ { 2 \kappa - 2 } +
	|  y | ^ { 2 \kappa - 2 } 
	\big) 
	|  x - y | ^ { 2 }.
\end{equation}
In addition, 
assume the initial data $ X _ 0 $ satisfies
\begin{equation}
 \label{eq:X_0-bound}
  \E \big[
   | X _ 0 | ^{ 2p ^ { * }  } 
	\big]
	 < 
	  \infty.
\end{equation}
\end{assumption}
According to \eqref{eq:SDEs-xy-contractive}, for any $q\in (0,\lfloor p ^ { * }\rfloor ]$,
there exists a sufficient small coefficient  
$ 0 < \epsilon _ 1 < \alpha_1 $
such that
\begin{equation}
 \label{eq:SDEs-x-contractive}
  \left \langle
	x ,  f ( x ) 
	 \right\rangle +
	  \tfrac{ 2 q- 1 } { 2 } 
	   \left\|
	    g ( x ) 
	     \right\| ^ { 2 } 
	\leq
	 - ( \alpha _ 1 - \epsilon _ 1 )
	   |  x  | ^ { 2 } +
	    c _ 0,	
\end{equation}    
where 
$ 
c _ 0 
=
\tfrac{ | f (0) | ^ 2 }{ 4 \epsilon _ 1 } +
 \tfrac{ ( 2p ^ {*} - 1 ) ^ 2 
  \left\| g (0) \right\| ^ 2 }{ 4 (p ^ { * }-q) } -
   \tfrac{ 2p ^ {*} - 1 }{2}
    \left\| g (0) \right\| ^ 2  $.
Furthermore, 
\eqref{eq:f-xy-condition} immediately implies
\begin{equation}
  \label{eq:f-x-condition}
   |  f ( x )  | ^ { 2 } 
	\leq
	  c _ 2
	   |  x  | ^ { 2 \kappa } +
	    c _ 3  ,
\end{equation}	
where  $ c _ 2 = \frac{ 2c _ 1 ( \kappa + 1 ) }{ \kappa }$ and 
$ c _ 3  = 2 | f ( 0 ) | ^ { 2 } +
\tfrac { 2c _1 ( \kappa - 1 ) }{ \kappa }$.

To show
%that Assumption \ref{ass:SDEs-xy-contractive} 
%is sufficient to 
the $2p$-th  moment ($p \geq 1$) boundedness  of $X(t)$,
 we first  establish the following lemma, 
which is a slight modification of
\cite [Lemma 8.1] {it1964stationary}.
\begin{lem}
 \label{lem:contractive-bound-moments}
Let 
$ r : [0,\infty) \rightarrow [0,\infty) $
and 
$ \xi : [0,\infty) \rightarrow [0,\infty) $ 
be nonnegative continuous functions.
If there exists a positive constant $\beta$ such that
 \begin{equation}
  \label{eq;r(t)-relation}
   r (t) -
 	r(s)	
 	 \leq
 	  - \beta
 	   \int _{ s } ^ {t}
 	    r (u)
 	     \, \dd u +
 	      \int _{ s } ^ {t}
 	       \xi (u)
 	        \, \dd u 
 \end{equation}
for any $ 0 \leq s < t < \infty $,
 then
 \begin{equation}
  \label{eq:r(t)-boundness}
 	r (t) 
 	\leq
 	 r (0) +
 	  \int _{ 0 } ^ {t}
 	   e  ^ { - \beta ( t - u )}
 	    \xi (u)
 	     \, \dd u.
 \end{equation} 
\end{lem}
\begin{comment}
    Write 
$
v _ 1 (t)
=
r(0) +
\int _{ 0 } ^ {t}
e  ^ { - \beta ( t - u )}
\xi (u)
\, \dd u $ .
It is easy to obtain that
\begin{equation}
 \begin{aligned}
  v _ 1 ^ {'} (t) 
   &
   =
    -\beta
     \int _{ 0 } ^ {t}
      e  ^ { - \beta ( t - u )}
       \xi (u)
        \, \dd u +
         \xi (t)
   \\
    &
    =
     -\beta
      ( v _ 1 (t) -
       r (0) ) +
        \xi (t).
 \end{aligned}
\end{equation}
So we have 
\begin{equation}
 \label{eq:v_1(t)-relation}
  v _ 1 (t) -
   v _ 1 (s)
  =
   -\beta
    \int _{ s } ^ {t}
     v _ 1 (u) -
      r (0)
       \, \dd u +
        \int _{ s } ^ {t}
         \xi (u)
          \, \dd u. 
\end{equation}
Now set 
$ 
v _ 2 (t) 
=
r(t) -
v _ 1 (t)$, 
according to \eqref{eq;r(t)-relation} 
and \eqref{eq:v_1(t)-relation}
we have that
\begin{equation}
 \begin{aligned}
  \label{eq:v_2(t)-relation}
   v _ 2 (t) -
    v _ 2 (s)
   &
   =
    r(t) -
     r(s) -
      ( v _ 1 (t) -
       v _ 1 (s) )
   \\
    &
     \leq
      -\beta
       \int _{ s } ^ {t}
        r (u) -
         v _ 1 (u) +
          r (0)
           \, \dd u 
   \\
    &
    =
     -\beta
      \int _{ s } ^ {t}
       v _ 2 (u) 
        \, \dd u .    
 \end{aligned}
\end{equation}
To prove \eqref{eq:r(t)-boundness}, 
it suffices to prove $ v _ 2 (t) \leq 0 $. 
If $ v _ 2 (t) > 0 $, 
for some $t$, 
then $ v _ 2 (0) = 0 $ 
implies that there exists an interval 
$ [s_1, t_1] \subset [0,\infty)$ 
with $ v _ 2 (t_1) > v _ 2 (s_1) $ 
and $  v _ 2 (t) > 0 $ on $ [s_1, t_1] $, 
which contradict \eqref{eq:v_2(t)-relation}. 
\end{comment}

Equipped with Assumption \ref{ass:SDEs-xy-contractive} 
and Lemma \ref{lem:contractive-bound-moments}, 
we shall give the moment boundedness of the  solution over infinite time, which is an interesting result independently.
\begin{thm}
 \label{thm:bound-moments}
Suppose Assumption \ref{ass:SDEs-xy-contractive}
hold. 
There exists a positive constant 
$ C _ 1 := C _ 1 ( p, \alpha _ 1, \epsilon_1, c _ 0 ) $ 
independent of $t$, 
such that
\begin{equation}
 \label{eq:X_t-moment-bound}
  \E \big[
    | X
	 ( t  ) | ^ { 2p }
	  \big]
	\leq
	 C _1
	  \E \big[
	   \big( 
	     1 +
	      | X _ 0 | ^ { 2p } 
	       \big)
	        \big],
	\qquad
	 1
	  \leq
	   p
	   <
	     p ^ { * }.
\end{equation}
%where $ X _ 0 $ is the initial value when
%$ t = 0 $.
\end{thm}
\textit{Proof of Theorem \ref{thm:bound-moments}}.	
By the $ \rm It \hat{o}$
formula,
one sees 
\begin{equation}
  \begin{aligned}
   \big(
    1 + | X(t) | ^ 2 
     \big) ^ p
	&
	=
	 \big( 
	  1 + | X_0 | ^ 2
	   \big) ^ p +
	    2p \int _ {0} ^ {t}
	     \big(
	      1 + | X(s) | ^ 2 
	       \big) ^ { p - 1 }
		    \left\langle
		     X(s),
		      f( X(s) )
		       \right \rangle
		        \, \dd s 
   \\
    &
	 \quad 
	  +
	   2p \int _ {0} ^ {t}
		\big(
		 1 + | X(s) | ^ 2 
		  \big) ^ { p - 1 }
		   \big\langle
		    X(s),
		     g( X(s) )
		      \, \dd W (s)
		      \big  \rangle 
	\\
	 &
	  \quad 
	   +	   
	    p  \int _ {0} ^ {t}
		 \big(
		  1 + | X(s) | ^ 2 
		   \big) ^ { p - 1 }
	 	    \left\|
		     g (X(s) ) 
		      \right\| ^ {2}
		       \, \dd s 
	\\
     &
	  \quad
	   +
	    2 p ( p - 1 )  
		 \int _ {0} ^ {t}
		  \big(
		   1 + | X(s) | ^ 2 
		    \big) ^ { p - 2 }
		     | X (s) ^ T
		      g (X(s) ) | ^ {2}
		       \, \dd s, 
	\end{aligned}
\end{equation}
which straightforwardly implies
\begin{equation}
 \label{eq:2p-power-bound}
   \begin{aligned}
	\big(
	 1 + | X(t) | ^ 2 
	 \big) ^ p
	&
	\leq
	 \big(
	  1 + | X_0 | ^ 2 
	   \big) ^ p  +
	    2p \int _ {0} ^ {t}
	     \big(
	      1 + | X(s) | ^ 2 
	       \big) ^ { p - 1 }
		    \big\langle
		     X(s),
		      f( X(s) )
		       \big \rangle
		        \, \dd s
	\\
	 &
	 \quad 
	  +
	   2p \int _ {0} ^ {t}
	    \big(
	     1 + | X(s) | ^ 2 
	      \big) ^ { p - 1 }
	       \left\langle
	        X(s),
	         g( X(s) )
	          \, \dd W (s)
	           \right \rangle
	\\
	 &
	 \quad 
	  +
	   p ( 2p - 1 )
		\int _ {0} ^ {t}
		 \big(
		  1 + | X(s) | ^ 2 
		   \big) ^ { p - 1 }
		    \left\|
		     g (X(s) ) 
		      \right\| ^ {2}
		       \, \dd s. 
	\end{aligned}
\end{equation}
According to \eqref{eq:SDEs-x-contractive}, 
we have
\begin{equation}
 \begin{aligned}
   \big(
    1 + | X(t) | ^ 2 
     \big) ^ p
   &
    \leq
	 \big(
	  1 + | X_0 | ^ 2 
	   \big) ^ p  -
	    2 p  ( \alpha_1 - \epsilon_1 )
	     \int _ {0} ^ {t}
	      \big(
	       1 + | X(s) | ^ 2 
	        \big) ^ p
	         \, \dd s 
   \\
    &
     \quad     
	  +
	   2 p(  c _ 0 + \alpha_1 - \epsilon_1 )
	    \int _{ 0} ^ {t}
	     \big(
	      1 + | X(s) | ^ 2 
	       \big) ^ {p-1}
	        \, \dd s
   \\
    &
     \quad 
      +
       2p \int _ {0} ^ {t}
        \big(
         1 + | X(s) | ^ 2 
          \big) ^ { p - 1 }
           \left\langle
            X(s),
             g( X(s) )
              \, \dd W (s)
               \right \rangle.	      
 \end{aligned}
\end{equation}
Then by 
the Young inequality, we obtain
\begin{equation}
 \begin{aligned}
  \label{eq:X_t-relation}
  \big(
    1 + | X(t) | ^ 2 
     \big) ^ p
   &
    \leq
    \big( 
	   1 + | X_0 | ^ 2 
	    \big) ^ p  -
	     p  ( \alpha_1 - \epsilon_1)
	      \int _ {0} ^ {t}
	       \big(
	        1 + | X(s) | ^ 2 
	         \big) ^ p
	          \, \dd s 
   \\
    &
     \quad       
	  +
	   \int _{ 0} ^ {t}
	    \frac{( 2 (  c _ 0 + \alpha_1 - \epsilon_1 ) ) ^ p}
	     {(   \alpha_1 - \epsilon_1  ) ^
	      { p - 1 } }
	       \, \dd s
   \\
	&
	 \quad 
	  +
	   2p \int _ {0} ^ {t}
	   \big(
	    1 + | X(s) | ^ 2 
	     \big) ^ { p - 1 }
	      \left\langle
	       X(s),
	        g( X(s) )
	         \, \dd W (s)
	          \right \rangle.          
 \end{aligned}
\end{equation}
For every integer $ n \geq 1 $, 
define the stopping time:
\begin{equation}
 \tau _ n 
 := 
  \inf 
   \{ 
    s \in [0, t]: 
     | X(s) |
      \geq n 
       \}.
\end{equation}
Clearly, 
$ \tau_n \uparrow t $ a.s..
Moreover, 
it follows from \eqref{eq:X_t-relation} 
and the property of the $ \rm It \hat{o}$ integral that
\begin{equation}
 \begin{aligned}
  & \E \Big[
   \big(
    1 +
     | X ( t \wedge \tau_n ) | ^ 2
      \big) ^ p 
       \Big] +
        p  ( \alpha_1 - \epsilon_1 )
         \E \bigg[
          \int _ {0} ^ {t \wedge \tau_n}
           \big(
            1 + | X(s) | ^ 2 
             \big) ^ p
              \, \dd s 
               \bigg]
   \\    
    &
     \leq  
      \E \Big[
       \big(
        1 +
        | X _ 0 | ^ 2
         \big) ^ p 
          \Big]  +
           \int _{ 0} ^ {t }
             \frac{( 2 (  c _ 0 + \alpha_1 - \epsilon_1 ) ) ^ p}
              {(   \alpha_1 - \epsilon_1 ) ^
               { p - 1 } }
                \, \dd s.
 \end{aligned}
\end{equation}
Letting $ n \to \infty $ 
and by the Fatou lemma, 
we get
\begin{equation}
 \begin{aligned}
  &
   \E \Big[
    \big(
     1 +
      | X ( t ) | ^ 2
       \big) ^ p 
        \Big] +
         p  ( \alpha_1 - \epsilon_1 )
         \int _ {0} ^ {t}
          \E \Big[
           \big(
            1 + | X(s) | ^ 2 
             \big) ^ p
              \Big]
               \, \dd s 
   \\            
    &
     \leq
      \E \Big[
       \big(
        1 +
         | X _ 0 | ^ 2
          \big) ^ p 
           \Big] +
            \int _{ 0} ^ {t }
             \frac{( 2 (  c _ 0 + \alpha_1 - \epsilon_1 ) ) ^ p}
               {(   \alpha_1 - \epsilon_1  ) ^
                { p - 1 } }
                 \, \dd s .
 \end{aligned}
\end{equation}
Taking
$r(t)=  \E \Big[
    \big(
     1 +
      | X ( t ) | ^ 2
       \big) ^ p 
        \Big]$,
        $ \beta = p  ( \alpha_1 - \epsilon_1 )  $ and 
$ \xi =  \frac{( 2 (  c _ 0 + \alpha_1 - \epsilon_1 ) ) ^ p}
          {(   \alpha_1 - \epsilon_1  ) ^
	       { p - 1 } }$
in  Lemma \ref{lem:contractive-bound-moments} shows
\begin{equation}
 \begin{aligned}
  \E \Big[
   \big(
     1 +
      | X ( t ) | ^ 2
       \big) ^ p 
        \Big]
   &
     \leq
      \E \Big[
       \big(
        1 +
         | X _ 0 | ^ 2
          \big) ^ p 
           \Big] +
	        \int _ {0} ^ {t}
	         e ^ {
		      - p ( \alpha_1 - \epsilon_1 )
		      ( t - u ) }
		       \frac{( 2 (  c _ 0 + \alpha_1 - \epsilon_1 ) ) ^ p}
		        {(   \alpha_1 - \epsilon_1  ) ^
		    	 { p - 1 } }
		          \, \dd u	
	\\
	 &
	  \leq
	   \E \Big[
	    \big(
	     1 +
	      | X _ 0 | ^ 2
	       \big) ^ p 
	        \Big] 	 +
	          \frac{( 2 (  c _ 0 + \alpha_1 - \epsilon_1 ) ) ^ p}
	           {p(   \alpha_1 - \epsilon_1  ) ^
	     	    p }
	             ( 1 -  e ^ {
	           	  - p  ( \alpha_1 - \epsilon_1 ) t })
   \\
	&
	\leq
	 C _ 1
	  \E \big[
	   \big( 
	    1 +
		 | X _ 0 | ^ {2p}
		  \big)
		   \big].
	\end{aligned}      
\end{equation}
Thus the proof is finished. 
$\square$

To establish a strong convergence theorem over an infinite time, 
we also rely on the following result.
\begin{thm}
	\label{lem:error-one-step-xy}
Let Assumption \ref{ass:SDEs-xy-contractive} 
hold.
For the representation
\begin{equation}
	\label{eq:X-xy-form}
	 X ( t, x; t + \theta ) -
	  X ( t, y; t + \theta)
  = 
   x - y
  +
    \mathcal{R} _ {t, x, y}
     (t + \theta ),
\end{equation}
we have for 
 $ 1 \leq p <\frac{   p ^ {*}    }{ \kappa } $,
\begin{equation}
	\label{eq:X-xy-2p}
	 \E \big[
	  | X ( t, x; t + h ) - X ( t, y; t + h ) | ^ { 2 p }
	   \big]
  \leq
	| x - y | ^{ 2p }
	 \exp( -2p \alpha _ 1 h),
\end{equation}	
\begin{equation}
 \label{eq:R-xy-2p}
  \E \big[ 
    |\mathcal{R} _ {t, x, y} (t + h ) | ^ { 2 p } \big]
  \leq
    C_2
	 \big(
	  1 +
	   | x | ^ { 2 \kappa - 2 } +
	    | y | ^ { 2 \kappa - 2 } 
	     \big) ^
	      { \frac{ p } { 2 } }
  	       | x - y | ^{ 2p } h ^ p,
\end{equation}	
where $ C _ 2 > 0 $ and do not depend on $ t, h$.
\end{thm}
\textit{Proof of Theorem
	\ref{lem:error-one-step-xy}.}
For $ s \geq t$,
we  introduce the process
\begin{center}
    $ \mathcal{D} _ {t, x, y}( s ) 
  := 
  X ( t, x; s) - X ( t, y; s) $
\end{center} 
and note that 
\begin{center}
  $ \mathcal{R} _ {t, x, y}  ( s ) 
  := 
  \mathcal{D} _ {t, x, y}( s )  - ( x - y ) $.  
\end{center}
Clearly,
\begin{equation}
    \dd \mathcal{D} _ {t, x, y}( s )
   =
    \big(
     f ( X ( t, x; s) ) - f( X( t, y; s) ) 
      \big)
       \, \dd s
    +
     \big(
       g ( X( t, x; s) ) - g( X ( t, y; s) )
      	\big)
      	 \, \dd W ( s ).
\end{equation}     	
%For the inequality \eqref{eq:X-xy-2p},
By applying the $ \rm It \hat{ o } $
formula,  
we obtain the following
for any $ t \leq s \leq t + h  $, 
\begin{equation}
 \begin{aligned}
	&
 \quad
 e ^ { 2 
	 \alpha _ 1  
	  p ( s - t ) }
	   | \mathcal{D} _ {t, x, y}( s ) | ^ { 2p } 
    -
    | \mathcal{D} _ {t, x, y}( t ) | ^ { 2p }
    \\
  &
  \leq
   \int  _ { t } ^ { s }
    2
     \alpha _ 1 p
	   e ^ { 2 
	   	\alpha _ 1 
	   	 p ( \theta - t ) }
		  | \mathcal{D} _ {t, x, y} ( \theta ) | ^{ 2p }
		   \, \dd \theta
  \\
   &
   \quad
    +
	 \int  _ { t  } ^ { s } 2p
	  e ^ { 2 
	   \alpha _ 1 
	     p ( \theta - t ) }
	      | \mathcal{D} _ {t, x, y}( \theta ) | ^{ 2p - 2 }
		   \left\langle
		    \mathcal{D} _ {t, x, y}( \theta ),
		     f ( X ( t, x; \theta ) ) -
		      f ( X ( t, y; \theta )  )
		       \right \rangle
		        \, \dd \theta
  \\
   &
   \quad
	+
	 \frac{ 2p (2p - 1 )}{ 2 }
	  \int  _ { t  } ^ { s }
		  e ^ { 2 
		   \alpha _ 1 
		    p ( \theta - t ) }
		     | \mathcal{D} _ {t, x, y} ( \theta ) | ^{ 2p - 2 }
		      \left\|
		       g ( X ( t, x; \theta ) ) - 
		        g ( X ( t, y; \theta ) ) 
		        \right\|^{2}
		         \, \dd \theta
  \\	
   &
   \quad
    +
	\int  _ { t  } ^ { s } 2p
	 e ^ { 2 
	 \alpha _ 1 
	  p ( \theta - t ) }
	   | \mathcal{D} _ {t, x, y} ( \theta ) | ^{ 2p - 2 }
		\left\langle
		 \mathcal{D} _ {t, x, y}( \theta ),
	      \big(
	       g ( X ( t, x; \theta ) ) - 
	        g ( X ( t, y; \theta ) )
	         \big)
	          \, \dd W ( \theta )
	           \right \rangle.  
	\end{aligned}
\end{equation}
Further, utilizing \eqref{eq:SDEs-xy-contractive} 
and taking expectations,
we have
\begin{equation}
	\E \big[
	  e ^ { 2 
	  \alpha _ 1 
	   p ( s - t ) }
	    | \mathcal{D} _ {t, x, y} ( s ) | ^ { 2p }
	     \big]
  \leq
   \E \big[
    | \mathcal{D} _ {t, x, y} ( t ) | ^ { 2p }
     \big].
\end{equation}
The estimate \eqref{eq:X-xy-2p} is straightforward.
Now we prove  \eqref{eq:R-xy-2p}.
By using the $ \rm It \hat{o} $ formula and 
\eqref{eq:SDEs-xy-contractive},
we obtain for $ \theta \geq 0 $,
\begin{equation}
 \begin{aligned}
   | \mathcal{R} _ {t, x, y}( t + \theta ) | ^{ 2p } 
    &
  	 \leq
  	   2p \int  _ { t } ^ { t + \theta }
  	    | \mathcal{R} _ {t, x, y}(s) | ^ { 2p-2 }
  	     \left\langle
  	       \mathcal{R} _ {t, x, y} ( s ),
  	        f ( X ( t, x; s ) ) - 
  	         f ( X ( t, y; s ) )
  	          \right \rangle
  	           \, \dd s 
  	  \\
  	   &
  	   \quad
  	    + 
  	     \frac{ 2p ( 2p-1 ) } { 2 }
  	      \int  _ { t } ^ { t + \theta }
  	       | \mathcal{R} _ {t, x, y}(s) | ^ { 2p-2 }
  	        \left\|
  	          g ( X( t, x; s) )-g ( X( t, y; s) )
  	           \right\| ^ { 2 }
  	            \, \dd s
   \\
  	&
  	\quad
  	 +
  	  2p \int  _ { t } ^ { t + \theta }
     | \mathcal{R} _ {t, x, y}(s) | ^ { 2p-2 }
  	   \left\langle
  	     \mathcal{R} _ {t, x, y} ( s ),
  	      \big( 
  	       g ( X ( t, x; s ) ) - 
  	        g ( X ( t, y; s ) ) 
  	         \big)
  	          \, \dd W ( s )
  	           \right \rangle
  	\\
  	 &
  	  \leq
  	   2p
  	    \int  _ { t } ^ { t + \theta }
  	      | \mathcal{R} _ {t, x, y}(s) | ^ { 2p-2 }
  	       \Big(
  	        \left\langle
  	         \mathcal{D} _ {t, x, y} ( s ),
  	          f ( X ( t, x; s ) ) - 
  	           f ( X ( t, y; s ) )
  	            \right \rangle
   \\
    &
   \quad
  	 + 
  	  \frac{  2p-1  } { 2 }
  	     \left\|
  	       g ( X( t, x; s) ) - 
  	        g ( X( t, y; s) )
  	         \right\| ^ { 2 }
  	          \Big)
  	           \, \dd s
    \\
  	 &
  	 \quad
  	   -2p
  	    \int  _ { t } ^ { t + \theta }
  	     | \mathcal{R} _ {t, x, y}(s) | ^ { 2p-2 }
  	      \left\langle
  	        x - y,
  	         f ( X ( t, x; s ) ) - f ( X ( t, y; s ) )
  	          \right \rangle 
  	           \, \dd s 
  	 \\
  	  &
  	  \quad
  	   +
  	    2p \int  _ { t } ^ { t + \theta }
       | \mathcal{R} _ {t, x, y}(s) | ^ { 2p-2 }
  	     \left\langle
  	     \mathcal{R} _ {t, x, y} ( s ),
  	       \big( 
  	        g ( X ( t, x; s ) ) - 
  	         g ( X ( t, y; s ) ) 
  	          \big)
  	           \, \dd W ( s )
  	            \right \rangle                        
   \\
     &
     \leq
      2p
       \int  _ { t } ^ { t + \theta }
         | \mathcal{R} _ {t, x, y}(s) | ^ { 2p-2 }
          \big| 
           f ( X ( t, x; s ) ) - 
            f ( X ( t, y; s ) ) 
             \big|
              | x - y  |
               \, \dd s
   \\
    &
    \quad
     +
      2p \int  _ { t } ^ { t + \theta }
      | \mathcal{R} _ {t, x, y}(s) | ^ { 2p-2 }
       \left\langle
        \mathcal{R} _ {t, x, y} ( s ),
          \big(
           g ( X ( t, x; s ) ) - 
            g ( X ( t, y; s ) ) 
             \big)
              \, \dd W ( s )
               \right \rangle.        
 	\end{aligned}
\end{equation}
According to the Cauchy-Bunyakovsy-Schwarz inequality (twice), 
\eqref{eq:f-xy-condition}, \eqref{eq:X_t-moment-bound} and 
\eqref{eq:X-xy-2p}, we get
\begin{equation}
  \begin{aligned}
  &
  \quad
   \E \big[
    | \mathcal{R} _ {t, x, y} ( t + \theta ) | ^ { 2p } 
     \big]
     \\
   &
  	  \leq
  	   2p | x - y |
  	    \int _ { t } ^ { t + \theta }
  	          \Big(
  	           \E \big[
  	            \big| 
  	             f ( X ( t, x; s) ) -
  	              f ( X ( t, y; s) ) 
  	               \big| ^ { p }
  	                \big]
  	                 \Big) ^ {\frac{ 2 }{ 2p } }
  	                  \Big(
  	                  \E \big[
  	                  | \mathcal{R} _ {t, x, y} ( s ) | ^ { 2p }
  	                  \big] 
  	                  \Big) ^ {  \frac{ 2p-2 } { 2p } }
  	                  \, \dd s  
  \\
   & 
   \leq
    2p
     \sqrt{ c _ 1 }
      | x - y |
       \int _ { t } ^ { t + \theta }
             \Big(
              \E \Big[
               \big(
                1 +
                 |  X ( t, x; s) | ^ { 2 \kappa - 2 }  +
                  |  X ( t, y; s)| ^ { 2 \kappa - 2 }
                   \big ) ^ { \frac { p } { 2 } }
  \\
   &
   \quad
    \big| 
     X ( t, x; s) - 
      X ( t, y; s) 
       \big| ^ { p }
        \Big] 
         \Big)^ { \frac{ 1 } { p } }
         \Big( 
         \E \big[
         | \mathcal{R} _ {t, x, y}( s ) | ^ { 2p }
         \big] 
         \Big) ^ { 1 - \frac{ 1 } { p } } 
          \, \dd s  
  \\
    & 
    \leq
     2p
      \sqrt{ c _ 1 }
       | x - y |
       \int _ { t } ^ { t + \theta }
            \Big(
             \E \Big[
              \big(
               1 +
                |  X ( t, x; s) | ^ { 2 \kappa - 2 } +
                 |  X ( t, y; s)| ^ { 2 \kappa - 2 }
                  \big ) ^ { p } 
                   \Big] 
                    \Big) ^ {  \frac{ 1 } { 2 p } } 
   \\
    &
    \quad
     \Big( 
      \E \big[ 
       \big| 
        X ( t, x; s) - 
         X ( t, y; s) 
          \big|  ^ { 2 p }
           \big]
            \Big) ^ { \frac{ 1 } { 2p } }
             \Big(
             \E \big[
             | \mathcal{R} _ {t, x, y} ( s ) | ^ { 2p }
             \big]
             \Big) ^ { 1 - \frac{ 1 } { p } } 
             \, \dd s  
   \\
  	&
  	\leq
  	  C ^ { ' }
  	   | x - y | ^ { 2 }
  	    \big( 
  	     1 +
  	      | x | ^ { 2 \kappa - 2 } +
  	       | y | ^ { 2 \kappa - 2 } 
  	       \big) ^ 
  	        { \frac{ 1 } { 2 } }
  	        \int _ { t } ^ { t + \theta }
  	         \Big(
  	          \E \big[ 
  	           | \mathcal{R} _ {t, x, y}( s ) | ^ { 2p } 
  	            \big]
  	             \Big) ^ { 1 - \frac{ 1 }{ p } }
  	              \, \dd s .
  \end{aligned}
\end{equation}
By employing the Gronwall inequality 
\cite[p.46, Theorem 8.3]{mao2007stochastic}, 
we obtain	
\begin{equation}
	\E \big[
	  | \mathcal{R} ( t + \theta ) | ^ { 2p } \big]
  \leq
	 C_2
	  \big(
	    1 +
	     | x | ^ { 2 \kappa - 2 } +
	      | y | ^ { 2 \kappa -2 } 
	       \big) ^ { \frac{ p }{ 2 } }
	        | x - y | ^ { 2p }
	         \theta ^ { p }.
\end{equation} 
The  proof of Lemma \ref{lem:error-one-step-xy} is thus completed.
$\square$

Now we present the long-time fundamental strong convergence theorem as follows.
\begin{thm}
	\label{thm:global-error}
$(\textit{A long-time fundamental strong convergence theorem})$ Suppose
	\begin{enumerate}
	\item[ ${ \textit( H _ 1 ) }$ ]
Assumption \ref{ass:SDEs-xy-contractive}  hold.
	\item[${ \textit( H _ 2 ) }$ ] 
The one-step approximation
 $Z ( t _ 0, X _ 0; h )$
is given by \eqref{eq:Z_t-def}
has the following orders of accuracy:
for some 
$ p \geq 1 $
there are 
$\eta _ 1 
    \geq 
      1  $, 
$ \eta _ 2 
    \geq 
     1 $
and $ 0 \leq h _ 0 \leq 1 $,
 $ q _ 1 > 1 $,
 $ \tfrac{ 1 } { 2 } 
    < 
     q _ 2 
     \leq 
      q _ 1 - 
       \tfrac{ 1 }{ 2 }$ 
such that for $ 0 < h\leq h _ 0 $,
the numerical method has,
respectively,
local weak and strong errors of order $ q _ 1 $ and $q_2$,
defined as
\begin{equation}
	\label{eq:one-step-weak-error}
	 \big|
	  \E \big[
	   X ( t _ 0, X _ 0; h ) - Z ( t _ 0, X _ 0; h)
	    \big]
	     \big|
  \leq
	 C_3
	  \E \big[
	   \big( 
	    1 +
	     | X _ 0 | ^ { \eta _ 1 } 
	      \big)
	       \big]
	        h ^ { q _ 1 },
\end{equation}
\begin{equation}
 \label{eq:one-step-strong-error}
	\Big(
	 \E \big[
	  | X ( t _ 0, X _ 0; h) - Z ( t _ 0, X _ 0; h)|  ^ { 2p } 
	   \big]
	    \Big) ^ { \frac{ 1 }{ 2p } }
  \leq
	 C _ 4
	 \Big(
	  \E \big[
	   \big(
	    1 + 
	     | X _ 0 | ^ { 2p \eta _ 2 } 
	      \big) 
	        \big] 
	         \Big) ^ 
	         { \frac{ 1 } { 2p } }
	        h ^ { q _ 2},
\end{equation}
where $ C _ 3, C _ 4 > 0$ and independent of $ h, t$.
\item[${\textit( H _ 3 ) }$]
The approximation $ Z _ k $
is given by
  \eqref{eq:Z_k-def} has finite moments,
i.e.,
for some $ p \geq 1 $
there exist 
$ \eta _ 3 
     \geq 
      1 $, 
$ h _ 0 > 0 $
and $ C _ 5 > 0 $ such that for all
 $ 0 < h\leq h _ 0 $ 
and 
   $ k \geq 0$,
\begin{equation}
 \label{eq:numerical-solution-bound}
  \E \big[
   | Z _ k | ^ { 2p }
    \big]
  \leq
   C_5
    \E \big[
     \big(
      1 +
       | X _ 0 | ^ { 2p \eta _ 3 } 
        \big)
         \big],
\end{equation}
where $ C _ 5 $ not depend on $ h, t, k$.
\end{enumerate}	
Then there exists a constant 
$\lambda:=
  \max
   \{ 
    \eta_1 \eta_3, 
     ( \tfrac{\kappa - 1 }{2} + 
      \eta_2) 
       \eta_3
        \} $,
for
   $ 
       h 
       \leq 
        h _ 1 
        :=
          \min
           \{ 
            \tfrac{ 1 }{ p \alpha _ 1 }, 
             h _ 0 
              \} $,
the global error is bounded as follows:
\begin{equation}
   \Big(
    \E \big[
	 |X_k  - Z _ k | ^ { 2p } 
	  \big]
	   \Big) ^ { \frac{ 1 } { 2p } }
  \leq
	 C 
  \Big(
    \E \big[
     \big(
      1 + 
      | X_0 | ^ { 2 \lambda p}
       \big)
         \big] 
          \Big)
           ^ {\frac{ 1 }{ 2p } }
           h ^ { q _ 2 - \frac{ 1 } { 2 } }, 
\end{equation}
where 
$ C 
  $ is a constant and do not depend on $h,t,k$.
\end{thm}
The proof is given in Appendix.

\section{Applications:
Strong Convergence Rate of the Backward Euler Method
over Infinite Time}
 \label{sec:backward} 
\noindent
In this section,
we will utilize the previously obtained strong convergence theorem to derive the strong convergence rate of the backward Euler scheme under  certain conditions.
The backward Euler method applied to SDEs \eqref{eq:form-SDE} 
takes the following form:
\begin{equation}
 \label{eq:discrete-backward-euler}
  Z _ k
 =
  Z _ {k-1}
 +
  f ( Z _ k ) h
 +
  g ( Z _ {k-1})
   \Delta W _ {k-1},
\end{equation}
where 
 $\Delta W _ {k-1}
   :=
      W ( t _ k )
   -  
      W ( t _ {k-1 }),
       k \geq 1  $.
Then the one-step approximation of 
\eqref{eq:discrete-backward-euler} 
reads:
\begin{equation}
 \label{eq:continuous-backward-euler}
	Z ( t, x; t + h )
   =
    x
   +
    \int _ {t} ^ { t + h }
     f ( Z ( t, x; t + h ) )
      \, \dd s
   +
    \int _ {t} ^ { t + h }
      g (x)
       \, \dd W(s).
\end{equation}
To apply the strong convergence theorem,
the first step is to ensure the $2p$-th  moment ($p \geq 1$) boundedness of $ Z_k$. 
%Hutzenthaler and Jentzen \cite{hutzenthaler2015numerical} have established that under the condition of  polynomial growth in the drift and diffusion coefficients of SDEs, the numerical solution remains bounded in the $p$-th  moment ($p \geq 1$) over a finite time. However,  ensuring the boundedness of the $2p$-th moment boundedness for $ Z_k $ over an infinite time horizon is challenging.
%
To this end, we make the following assumptions.
\begin{assumption}
 \label{ass:g-xy-condition}
Suppose that the diffusion coefficient of SDEs 
\eqref{eq:form-SDE} satisfies the global Lipschitz condition.
Namely,
there exists a constant $ \beta_1 > 0 $ such that
\begin{equation}
 \label{eq:g-xy-lipschitz}
  \left\|
    g(x) - g(y)
     \right\| ^ 2
   \leq
    \beta_1
     | x - y | ^ 2.
\end{equation}
\end{assumption}

According to Assumtions
\ref{ass:SDEs-xy-contractive},
\ref{ass:g-xy-condition},
there exist some positive constants 
$\alpha_2, 
  \beta _ 2 
   $ 
and 
 $ 2 
    \alpha_2 
     > 
      ( 2p ^ {*}  - 1 ) 
       \beta _ 2 $, 
such that
\begin{equation}
   \label{eq:SDEs-x-contractive-backward}
	\left \langle
	 x ,  f ( x ) 
	 \right \rangle  
   \leq
	-\alpha_2
	 |  x  | ^ { 2 } +
	  c _ 4,
\end{equation}
\begin{equation}
 \label{eq:g-x-lipschitz}
  \left\|
	g ( x )
	 \right\| ^ { 2 } 
   \leq
	\beta _  2 
	 |  x  | ^ { 2 } +
	  c _ 5,
\end{equation}
where $ c _ 4, c _ 5 > 0$ are independent of $x$.

\begin{lem}
 \label{lem:one-step-2p-relation}
If there exist some positive constants $ \mu$, 
$\nu$ such that $ 1 - \mu h > 0 $ 
and  the sequence $\{Y_k\} _ { k \in \N }$ satisfies
\begin{equation}
 \label{eq:one-step-2p-relation}
	| Y _ k | ^ { 2 p } 
	 \leq
	   ( 1 - \mu h )
	     | Y _{ k - 1 } | ^ { 2 p } +
	      \nu h  
	       | Y _{ k - 1 } | ^ { 2 p -2  }
\end{equation}
for any 
$ p \geq 1 $ 
and 
$ k \in \N $,
then we obtain
\begin{equation}
 \label{eq:Y-k-bound}	
   | Y _ k | ^ { 2 p } 
    \leq
     C \big(
      1 + | Y _ 0 | ^ { 2 p } 
       \big),     	
\end{equation} 
where $ C $ only depends on 
$ p, \mu, \nu$ 
and
 $ Y _ 0 $ 
 is the initial value.
\end{lem}
\textit{Proof of Lemma \ref{lem:one-step-2p-relation}. }
When $ p = 1 $, 
by 
\eqref{eq:one-step-2p-relation} 
we obtain
\begin{equation}
 \begin{aligned}
  \label{eq:one-step-2-result}
   | Y _ k | ^ 2
    &
     \leq
      ( 1 - \mu h ) ^ k
       | Y _ 0 | ^ 2 +
        \big(
         1 +
          ( 1 - \mu h ) +
           \cdots +
            ( 1 - \mu h ) ^ { k - 1 }
             \big)
              \nu h  
   \\
    &
     =
      ( 1 - \mu h ) ^ k
       | Y _ 0 | ^ 2 +
        \frac{ 1 - ( 1 - \mu h ) ^ k }{ \mu h }
         \nu h  
   \\
    &
     \leq
       C \big(
        1 + | Y _ 0 | ^ 2
        \big).       
 \end{aligned}
\end{equation}
For 
$ p > 1 $, 
we can utilize the Young inequality to obtain:
\begin{equation}
  \begin{aligned}
   | Y _ k | ^ { 2 p }
    &
     \leq
       \big( 1 - \frac{\mu}{2} h 
       \big)
       | Y _{ k - 1 } | ^ { 2 p } +
       \nu ^ p  
        \big( \frac{ 2 p - 2	}{ \mu p } 
        \big) 
        ^ { p - 1 }	h.
   \end{aligned}
\end{equation}
Following the same procedure as in
\eqref{eq:one-step-2-result}, 
one can straightforwardly derive the inequality
\eqref{eq:Y-k-bound}.
$\square$
\begin{thm}
 \label{thm:bound-backward-euler}
Let  Assumptions
  \ref{ass:SDEs-xy-contractive}
   and
 \ref{ass:g-xy-condition} 
 hold.
For 
$ p \in [ 1, \lfloor p ^ {*} \rfloor]$ 
and 
$ k \in \N $,
there exists a constant
    $ K := K ( \alpha_2, \beta _ 2, c _ 4, c _ 5, p ) > 0 $ 
such that
\begin{equation}
	\label{eq:bound-backward-euler}
 	 \E \big[ | 
 	   Z _ k
 	   |^ {2p} 
 	    \big]
   \leq
 	K
 	 \E \big[
 	  \big(
 	   1 +
        | X _ 0 | ^ {2p}
         \big)
          \big].
\end{equation}
\end{thm}
In the following proof and throughout the rest of the paper,
the symbol
  $ K $ 
will represent a generic constant that may vary from line to line.
Specifically, 
the factors upon which $ K $ depends may differ,
but it remains independent of $ t, h, k$.
     
 \textit{Proof of Theorem 
 	\ref{thm:bound-backward-euler}.} 
From
\eqref{eq:discrete-backward-euler}
 and \eqref{eq:SDEs-x-contractive-backward}, 
we have
\begin{equation}
  \label{eq:relate-Z_k-Z_{k-1}}
 	\begin{aligned}
 		( 1 + 2 \alpha_2 h ) ^ p
 		 | Z _  k| ^ {2p}
  \leq
   \big(
    | Z _ {k-1} +
     g ( Z _ { k - 1 } )
      \Delta W _ { k - 1 } | ^ 2
  +
   K h
   \big) ^ p.
 	\end{aligned}
\end{equation}  

Let us first consider the case $ p = 1 $. 
Taking  expectations of \eqref{eq:relate-Z_k-Z_{k-1}} 
and employing  \eqref{eq:g-x-lipschitz}, 
we get
\begin{equation}
 \label{eq:iteration-Z_k-Z_{k-1}}
  \begin{aligned}
	( 1 + 2 \alpha_2 h )
	 \E \big[
	  | Z _ k | ^ {2}
	   \big]
  \leq
	( 1 +
     \beta _ 2 h)
      \E \big[
        | Z _ { k - 1 } | ^ {2}
         \big]
  +
   K h.
  \end{aligned}
\end{equation} 
Since $ 2 \alpha_2 > \beta _ 2 $,
we have
 $\frac{ 1 + \beta _ 2 h }{ 1 + 2 \alpha_2 h } = 1 - \gamma h $,
where
$ \gamma = \frac{ 2 \alpha_2 - 
	\beta _ 2 }{1 + 2 \alpha_2 h } > 0 $.
By \eqref{eq:iteration-Z_k-Z_{k-1}} and Lemma 
\ref{lem:one-step-2p-relation},
we  derive that
\begin{equation}
	\begin{aligned}
	 \E \big[
	  | Z _ {k} | ^ {2}
	   \big]
  &
   \leq
    ( 1 -
   \gamma 
     h) 
     \E \big[
      | Z _ { k - 1 } | ^ {2}
       \big]  + 
         K h
  \\
   &
    \leq
	 K
	  \E \big[
	   \big( 
	    1 +
         | X _ 0 | ^ {2} 
          \big)
           \big].
	\end{aligned}
\end{equation}
For integer $ p > 1 $,
we note that:
\begin{equation}
\label{eq:form-I_1}
	\begin{aligned}
	 J_1
  &
  :=
   \E \big[
    \big(
     | Z _ { k - 1 } +
      g ( Z _ { k - 1 } )
       \Delta W _ { k - 1 } 
       | ^ 2
  +
   K h
   \big) ^ { p  }
    \big|
     \mathcal{F} _ { t _ {k-1} }
      \big]
  \\
   &
    \leq
     \E \big[
      | Z _ { k - 1 } +
       g ( Z _ { k - 1 } )
        \Delta W _ { k - 1 }  | ^ 
         { 2 p }
         \big|
          \mathcal{F} _ { t _ {k-1} }
           \big]
   +
    K h
    \big( 1 +
     | Z _ { k - 1 } | ^ { 2 ( p  - 1 ) }
      \big)
  \\
   &
   =
    \E \Big[
     \sum _ { l = 0 } ^ { 2 p }
      C _ {2p} ^ l
       | g ( Z _ { k - 1 } )
        \Delta W _ { k - 1 } | ^{ l }
         | Z _ { k - 1 } | ^
          { 2 p - l }
           \Big|
            \mathcal{F} _ { t _ {k-1} }
             \Big]
  +
   K h
   \big( 1 +
    | Z _ { k - 1 } | ^ { 2 ( p - 1 ) }
    \big) 
  \\
   & 
   \leq
    ( 1 
     + c (p,\beta_2,h) ) 
     | Z _ { k - 1 } | ^ { 2 p }
  +
   K h
   \big( 1 +
    | Z _ { k - 1 } | ^ { 2 ( p  - 1 ) }
     \big),
	\end{aligned}
\end{equation}       
where 
$
 c (p,\beta_2,h) := 
   \sum _ { l = 1 } ^ {  p }
    C _ {2p} ^ {2l} h ^ l \beta_2 ^l
     ( 2l - 1 )!!    
$.      
By taking the conditional expectations of 
\eqref{eq:iteration-Z_k-Z_{k-1}}
and then applying \eqref{eq:form-I_1},
we see
\begin{equation}
  \E \big[
   | Z _ k | ^ { 2 p  }
    \big]
   \leq
    (1-\gamma h )
     \E \big[
      | Z _ { k - 1 } | ^ 
       { 2 p }
       \big] +
         K h
          \E \big[
           \big( 1 +
            | Z _ { k - 1 } | ^ { 2 ( p- 1 ) }
             \big) 
              \big],
\end{equation}
where
$ \gamma
 =
    \frac{(1 +
	 2 \alpha_2 h) ^
	 { p } -
	  (1 +  c (p,\beta_2,h)) } 
      { h
       (1 +
	    2 \alpha_2 h) ^ 
	     { p  } }
           $,  
Verifying $ \gamma > 0 $ is straightforward.
Therefore, by applying Lemma
 \ref{lem:one-step-2p-relation}  and 
 Theorem \ref{thm:bound-backward-euler}, 
 we have established the result 
 for integer $p$. 
For non-integer $p$, 
we can derive the result using the Young inequality. 
Thus the proof is finished.
 $\square$       
 
In order to analyze the strong convergence rate,
we also introduce an auxiliary one-step approximation,
\begin{equation}
  \label{eq:auxiliary one-step approximation}
  	Z _ E( t, x; t + h )
  :=
    x
  +
    \int _ {t} ^ { t + h }
      f(x)
        \, \dd s 
  +
    \int  _ {t} ^ { t + h }
      g(x)
        \, \dd W ( s ), 
\end{equation}
which can be regarded as a one-step approximation of the Euler-Maruyama method.
Subtracting
       \eqref{eq:auxiliary one-step approximation}
from \eqref{eq:continuous-backward-euler} yields 
\begin{equation}
\label{eq:I_f-def}
  	Z ( t, x; t + h )
  :=
    Z _ E( t, x; t + h )
  +\underbrace{
     \int _ {t} ^ { t + h }
       f ( Z( t, x; t + h ) )
   -
     f (x)
      \, \dd s }_{:=I_f}.
\end{equation}
To derive the one-step error of the backward Euler scheme \eqref{eq:continuous-backward-euler},
we require several lemmas. 
\begin{lem}
	\label{lem:one-step-2p-bound}
Suppose Assumptions \ref{ass:SDEs-xy-contractive}
 and
\ref{ass:g-xy-condition} hold.
For any
$ 1 
  \leq 
   p 
    <
     \frac{ p ^ {*}  }{\kappa} $, 
one gets that
\begin{equation}
 \label{eq:one-step-2p-bound}
  \begin{aligned}
   \E \big[
	| X ( h ) - X _ 0 | ^ {2p}
	 \big]
   \leq
	 C _ 6
	  \E \big[
	   \big(
	    1 +
	     | X _ 0 | ^ { 2 \kappa p } 
	      \big)
	       \big]
		    h ^ {p},
		\end{aligned}
	\end{equation}	
where $ C _ 6 $ does not depend on 
$h,t$.
\end{lem}
\textit{Proof of Lemma 
	\ref{lem:one-step-2p-bound}.}
By applying
\eqref{eq:f-x-condition},
\eqref{eq:g-x-lipschitz} 
and
\eqref{eq:X_t-moment-bound},
and then combining the H{\"o}lder inequality with the moment inequality
\cite[Theorem 7.1]{mao2007stochastic},
we can obtain that
\begin{equation}
 \begin{aligned}
  \E \big[
   | X ( h ) - X _ 0 | ^ {2p}
	\big]
   &
   =
	\E \Big[
	 \Big|
	  \int _ {0} ^ {h}
	   f ( X (s) )
		\, \dd s +
		 \int _ {0} ^ {h}
		  g ( X (s) )
		   \, \dd  W(s)
		    \Big| ^ {2p}
		     \Big]
%   \\
%	&
%	\leq     
%	 C _ p
%	  \E \Big[
%	   \Big|
%		\int _ {0} ^ {h}
%		 f ( X (s) )
%		  \, \dd s  
%		   \Big|^ {2p}
%		    \Big] +
%		     C _ p
%		      \E \Big[
%		       \Big|
%		        \int _ {0} ^ {h}
%		         g ( X (s) )
%		          \, \dd  W(s)
%		           \Big| ^ {2p}
%		            \Big]
   \\
	&
	\leq 
	 C _ p 
	  h ^ { 2p - 1 }
	   \int _ {0} ^ {h}
		\E \big[
		 | f ( X (s) ) | ^ {2p}
		 \big]
		  \, \dd s +
		   C _ p
		    h ^ { p - 1 }
		     \int _ {0} ^ {h}
		      \E \big[
		       \left\|
		        g ( X (s) )
		         \right\|^ {2p}
		          \big]
		           \, \dd  s
   \\
	&
	\leq 
	 C _ p
	  h ^ { 2p - 1 }
	   \int _ {0} ^ {h}
		\E \Big[
		 \big(
		  c _ 3 +
		   c _ 2
		    | X (s) | ^ { 2 \kappa }
		     \big) ^ {p}
		      \Big]
		       \, \dd s +
		        C _ p
		         h ^ { p - 1 }
		          \int _ {0} ^ {h}
		           \E \Big[
		            \big(
		             c _ 5+
		              \beta _ 2
		               | X (s) | ^ 2
		                \big) ^ {p}
		                 \Big]
		                  \, \dd s
   \\
	& 
	\leq
	 C _ 6
	  \E \big[
	   \big(
	    1 +
		 | X _ 0 | ^ { 2 \kappa p } 
		  \big)
		    \big]
	  	    h ^ {p}.
	\end{aligned}
\end{equation}
Here
and in the following, 
the letter
$ C _ p $  is used to denote a generic positive constant, 
which only depend on $p$ 
and may vary for each appearance.
$\square$
\begin{lem}
	\label{lem:one-step-backward-euler-2p-bound}
Let Assumptions \ref{ass:SDEs-xy-contractive}
 and
\ref{ass:g-xy-condition} be satisfied.
For any 
$ 1 
   \leq 
    p 
     \leq
      \frac{ \lfloor p ^ {*} \rfloor }{\kappa} $,
one can see that
\begin{equation}
  \begin{aligned}
   \E \big[
	\big|
	  Z ( t, x; t + h ) - x
	  \big| ^ {2p}
	   \big]
	    \leq
	     C _ 7 
	      \big(
	       1 +
	        | x | ^ 
	         { 2 \kappa p } 
	          \big) 
	           h ^ {p},
  \end{aligned}
\end{equation}
where $ C _ 7 $ is independent of $h,t$.
\end{lem}      
\textit{Proof of Lemma \ref{lem:one-step-backward-euler-2p-bound}}.
Applying  \eqref{eq:f-x-condition},
 \eqref{eq:g-x-lipschitz},
 \eqref{eq:bound-backward-euler}, the  H{\"o}lder inequality 
and the moment inequality,
yields
\begin{equation}
 \begin{aligned}
  \E \big[
   \big|
    Z ( t, x; t + h ) - x
     \big| ^ {2p}
      \big]
  &
  =
   \E \Big[
    \Big|
     \int _ {t} ^ { t + h }
      f ( Z ( t, x; t + h ) )
       \, \dd s 
  +
   \int _ {t} ^ { t + h }
    g ( x )
     \, \dd W ( s )
      \Big| ^ {2p}
       \Big]
  \\
   &
    \leq   
      C _ p h ^ { 2 p - 1 }
       \int _ {t} ^ { t + h }
        \E \big[
         \big|
          f ( Z ( t, x; t + h ) )
           \big| ^ {2p}
            \big]
             \, \dd s
  +
   C _ p h ^ { p - 1 }
   \int _ {t} ^ { t + h }
    \E \big[
     \left\|
      g (x)
       \right\| ^ { 2 p }
        \big]
         \, \dd s
  \\
   &
   \leq
     C _ 7 
      \big( 1 +
       | x | ^ 
        { 2 \kappa p } 
         \big) 
          h ^ {p}.
	\end{aligned}
\end{equation} 
Thus the proof is finished. 
$\square$
\begin{lem}
	\label{lem:auxiliary-backward-euler-weak-error}
Let Assumptions 
\ref{ass:SDEs-xy-contractive} and
\ref{ass:g-xy-condition} 
be fulfilled.
Then
\begin{equation}
 \label{eq:auxiliary-backward-euler-weak-error}
  \begin{aligned}
   \big|
    \E [
     X ( t, x; t + h ) -
       Z _ E ( t, x; t + h ) ]
       \big|
  &
   \leq
     K 
     \big( 1 +
      | x | ^ { 2 \kappa - 1 } 
       \big)
        h ^ { \frac{3}{2} },
  \\
    \big|
     \E [
      X ( t, x; t + h ) -x ]
       \big|
   &
    \leq
     K 
      \big( 1 +
       | x | ^ { \kappa }
        \big)
          h,
   \\
     \big|
      \E [
       Z _ E ( t, x; t + h ) -x ]
        \big|
   &
    \leq
     K 
      \big( 1 +
       | x | ^ {\kappa} 
        \big) h. 
  \end{aligned}
\end{equation}
\end{lem}
\textit{Proof of Lemma \ref{lem:auxiliary-backward-euler-weak-error}}.
For the first item in 
\eqref{eq:auxiliary-backward-euler-weak-error},
we can utilize a variant of 
\eqref{eq:f-xy-condition},
   \eqref{eq:X_t-moment-bound}
and 
\eqref{eq:one-step-2p-bound}, 
as well as the H{\"o}lder inequality to acquire that
\begin{equation}
 \begin{aligned}
  \big|
   \E [
     X ( t, x; t + h ) -
       Z _ E ( t, x; t + h ) ]
      \big|
  &
  =
   \Big|
    \E \Big[
     \int _ {t} ^ { t + h }
      f ( X (s) ) - f (x)
       \, \dd s
  +
   \int _ {t} ^ { t + h }
    g ( X (s) ) - g (x)
     \, \dd W (s)
      \Big]
       \Big|
%  \\
%   &
%    \leq
%     \int _ {t} ^ { t + h }
%      \E \big[
%       \big|
%        f ( X(s) ) - f(x)
%         \big|
%          \big]
%           \, \dd s
  \\
   & 
    \leq
     \sqrt{ c _ 1 }
      \int _ {t} ^ { t + h }
       \E \Big[
        \big(
         1 +
          | X (s) | ^
           { 2 \kappa - 2 }
  +
   | x | ^
    { 2 \kappa - 2 }
     \big) ^
      {\frac{1}{2} }
       \big| 
        X (s) - 
         x 
          \big|
           \Big]
            \, \dd s
  \\
   & 
    \leq
     \sqrt { c _ 1 }
      \int _ {t} ^ { t + h }
       \Big(
        \E \Big[
         \big(
          1 +
           | X (s) | ^
            { 2 \kappa - 2 } +
             | x | ^
              { 2 \kappa - 2 } 
               \big) ^ 
                {\frac{ 2 \kappa - 1 }{ 2 \kappa - 2 }}
                 \Big] 
                  \Big) ^
                   { \frac{ \kappa - 1 }{ 2 \kappa - 1 } }
                   \\
                   &
                   \quad
                    \Big(
                     \E \Big[
                      \big|
                       X (s) - x 
                        \big| ^ 
                         {\frac{ 2 \kappa - 1 }{\kappa}}
                          \Big]
                           \Big) ^
                            { \frac{ \kappa }{ 2 \kappa - 1 } }
                             \, \dd s
  \\
   &
    \leq 
      K 
      \int _ {t} ^ { t + h }
       ( 1 +
        | x | ^
         { \kappa - 1 } )
          ( 1 +
           | x | ^
            {\kappa} )
             ( s - t ) ^
              { \frac{1}{2} }
               \, \dd s 
  \\
   & 
   \leq 
    K
     \big( 1 +
      | x | ^
      { 2 \kappa - 1 }
       \big)
        h ^ { \frac{3}{2} }.
 \end{aligned}
\end{equation}  
%Then,
%we bound  the second item in  \eqref{eq:auxiliary-backward-euler-weak-error},
%by applying the inequalities \eqref{eq:f-x-condition} and \eqref{eq:X_t-moment-bound},
%which leads to
%\begin{equation}
% \begin{aligned}
%  \big|
%   \E [
%    \rho _ { X, x } ]
%     \big|
%  &
%  =
%   \Big|
%    \E \Big[
%     \int _ {t} ^ { t + h }
%      f ( X(s) )
%       \, \dd s
%  +
%   \int _ {t} ^ { t + h }
%    g( X (s) )
%     \, \dd W (s)
%      \Big]
%       \Big|
%  \\
%   &
%    \leq
%     \int _ {t} ^ { t + h }
%      \E \big[
%       \big|
%        f ( X (s) ) 
%        \big|
%         \big]
%          \, \dd s
%  \\
%   & 
%   \leq
%    \int _ {t} ^ { t + h }
%     \E \big[
%      \big(
%       c _ 2 
%  +
%   c _ 3
%    | X (s) | ^
%     { 2 \kappa }
%      \big) ^
%       { \tfrac{1}{2} }
%        \big]
%         \, \dd s
%  \\
%   &
%   \leq
%    K 
%     \big( 1 +
%     | x | ^ {\kappa}
%      \big) h.
% \end{aligned}
%\end{equation}
%In addition, 
%using \eqref{eq:f-x-condition} and \eqref{eq:X_t-moment-bound}, 
%a similar result can be derived for
%  $ \big|
%     \E [
%      \rho _ { Z _ A, x } ]
%        \big|$.
%Thus the proof is accomplished.
Furthermore,
utilizing \eqref{eq:f-x-condition} 
and \eqref{eq:X_t-moment-bound}, 
the remaining terms in
\eqref{eq:auxiliary-backward-euler-weak-error} can be analogously validated.
Thus the proof is accomplished. 
$\square$
\begin{lem}
	\label{lem:one-step-backward-weak-error}
Under Assumptions
	 \ref{ass:SDEs-xy-contractive} and
	 \ref{ass:g-xy-condition}, 
the one-step backward Euler scheme \eqref{eq:continuous-backward-euler} obeys 
\begin{equation}
 \label{eq:one-step-backward-weak-error}
  \begin{aligned}
   \big|
    \E [
      X ( t, x; t + h ) -
       Z ( t, x; t + h ) ]
       \big|
   &
    \leq
     K 
      \big( 1 +
       | x | ^
        { 2\kappa - 1 }
         \big)
          h ^ { \frac{3}{2} },
  \\
   \big|
    \E 
     [ Z ( t, x; t + h )]
      \big|
  &
   \leq
    K 
     \big( 1 +
     | x | ^
      {\kappa} 
       \big) h .
  \end{aligned}
\end{equation}
\end{lem}
\textit{Proof of Lemma \ref{lem:one-step-backward-weak-error}}.
We will begin by discussing the following inequality:
\begin{equation}
 \label{eq:decomposition-one-step-backward-weak-error}
  \begin{aligned}
   \big|
    \E [
      X ( t, x; t + h ) -
       Z  ( t, x; t + h ) ]
       \big|
  &
  =
   |
    \E [
      X ( t, x; t + h ) -
       Z _ E ( t, x; t + h ) -
       I _ f ]
        |
  \\
   & 
   \leq
    \big|
     \E [
       X ( t, x; t + h ) -
       Z _ E ( t, x; t + h )]
        \big|
  +
  \big|
   \E 
    [ I _ f ]
     \big|.
  \end{aligned}
\end{equation}
The first term on the right-hand side of
 \eqref{eq:decomposition-one-step-backward-weak-error}
is estimated in  Lemma \ref{lem:auxiliary-backward-euler-weak-error}, 
yielding:
\begin{equation}
 \label{eq:weak-error-auxiliary-backward}
  \begin{aligned}
   \big|
    \E [
      X ( t, x; t + h ) -
       Z _ E ( t, x; t + h ) ]
       \big|
  &
  \leq
   K 
    \big( 1 +
     | x | ^
      { 2 \kappa - 1 } 
       \big)
        h ^ { \frac{3}{2} }.
  \end{aligned}
\end{equation}
For the second item on the right-hand side of
 \eqref{eq:decomposition-one-step-backward-weak-error},
by using
  \eqref{eq:f-xy-condition},
  along with
Theorem \ref{thm:bound-backward-euler}
and  Lemma \ref{lem:one-step-backward-euler-2p-bound},
as well as  applying the H{\"o}lder inequality,
one can obtain
\begin{equation}
  \begin{aligned}
   \big|
    \E 
     [ I _ f ]
      \big|
  &
  =
   \Big|
    \E \Big[
     \int _ {t} ^ { t + h }
      f ( Z ( t, x; t + h ) ) -
       f (x)
        \, \dd s
         \Big]
          \Big|
  \\
   &
    \leq
     \int _ {t} ^ { t + h }
      \E \big[
       \big|
        f ( Z ( t, x; t + h ) ) -
         f (x)
          \big|
           \big]
            \, \dd s 
  \\
   &
    \leq
     \sqrt { c _ 1 }
      \int _ {t} ^ { t + h }
       \E \Big[
        \big( 1 +
         |
          Z ( t, x; t + h )
           | ^
            { 2 \kappa - 2 }
  +
   | x |^
    { 2 \kappa - 2 }
     \big) ^
      { \frac{1}{2} }
       \big|
        Z ( t, x; t + h ) - x
         \big|
          \Big]
           \, \dd s
  \\
   &
    \leq
     \sqrt { c _ 1 }
      \int _ {t} ^ { t + h }
       \Big(
        \E \Big[
         \big(
          1 +
           |
            Z ( t, x; t + h )
             | ^
              { 2 \kappa - 2 } +
               | x | ^
                { 2 \kappa - 2 }
                 \big) ^ 
                  { \frac{2 \kappa - 1 }{2 \kappa - 2 }}
                   \Big] 
                    \Big) ^
                     { \frac{\kappa - 1}{2 \kappa - 1 } }
                     \\
                     &
                     \quad
                      \Big(
                       \E \Big[
                        | Z ( t, x; t + h ) -
                         x | ^ 
                          {\frac{2 \kappa - 1}{ \kappa }}
                           \Big]
                            \Big) ^
                             { \frac{\kappa}{2 \kappa - 1} }
                              \, \dd s
% \\
%  & 
%  \leq
%    K
%     \int _ {t} ^ { t + h }
%      \big( 
%       1 +
%        | x | ^
%        { \kappa - 1 } 
%         \big)
%          \big( 1 +
%           | x | ^
%            {\kappa} 
%             \big)
%              h ^ { \tfrac{1}{2} }
%              \, \dd s
 \\
  & 
  \leq
   K 
   \big( 1 +
    | x | ^
     { 2 \kappa - 1 } 
      \big)
       h ^ { \frac{3}{2} },
  \end{aligned}
\end{equation}
where $I_f$ is given by \eqref{eq:I_f-def}.
Combining this with \eqref{eq:weak-error-auxiliary-backward},
we obtain the first inequality in 
\eqref{eq:one-step-backward-weak-error}.
%Then, 
%we consider the second item in
% \eqref{eq:one-step-backward-weak-error},
%by utlizing \eqref{eq:f-x-condition}
%and Theorem \ref{thm:bound-backward-euler}, 
%which derives to 
%\begin{equation}
% \begin{aligned}
%  |
%   \E 
%    [ \rho _ { Z, x } ]
%     |
%  &
%  =
%   \Big|
%    \E \Big[
%     \int _ {t} ^ { t + h }
%      f ( Z ( t, x; t + h ) )
%       \, \dd s
%  +
%   \int _ {t} ^ { t + h }
%    g (x)
%     \, \dd W (s)
%      \Big]
%       \Big|
%  \\
%   &
%    \leq
%     \int _ {t} ^ { t + h }
%      \E \big[
%       |
%        f ( Z ( t, x; t + h ) )
%         |
%          \big]
%           \, \dd s
%  \\
%   & 
%    \leq
%     \int _ {t} ^ { t + h }
%      \E \Big[
%       \big(
%        c _ 2 + c _ 3
%         |
%          Z ( t, x; t + h )
%           | ^
%            { 2 \kappa }
%             \big) ^
%              { \tfrac{1}{2} }
%               \Big]
%                \, \dd s
%  \\
%   &
%    \leq
%     K 
%      \big(
%        1 +
%         | x | ^
%          {\kappa} 
%           \big) 
%            h.
% \end{aligned}
%\end{equation}
Thanks to \eqref{eq:f-x-condition}
and Theorem \ref{thm:bound-backward-euler}, 
the second item in
 \eqref{eq:one-step-backward-weak-error} 
becomes apparent.
 $\square$
\begin{lem}
 \label{lem:one-step-backward-strong-error}
Under Assumptions
	 \ref{ass:SDEs-xy-contractive},
	 \ref{ass:g-xy-condition},
for any 
$ 1 
   \leq
    p
     \leq 
      \frac{ 
      \lfloor p ^ {*} 
      \rfloor
      }{2\kappa - 1} $,
we obtain 
\begin{equation}
 \label{eq:one-step-backward-strong-error}
  \begin{aligned}
   \Big(
    \E \big[
     |
      I _ f 
       | ^ {2p}
        \big]
         \Big) ^
          { \frac{1}{2p} }
  &
  \leq
   K
   \Big(
    1 +
     | x | ^
      { ( 4 \kappa - 2 ) p }
       \Big) ^
        { \frac{1}{2p} }
         h ^ { \frac{3}{2} },
  \\
   \Big(
    \E \big[
     |
       X ( t, x; t + h ) -
       Z  ( t, x; t + h )
        | ^ {2p}
         \big]
          \Big) ^
           { \frac{1}{2p} }
  &
  \leq
   K
   \Big(
    1 +
     | x | ^
      { ( 4 \kappa - 2 ) p }
       \Big) ^
        { \frac{1}{2p} } h,
  \end{aligned}
\end{equation}	
where $I_f$ is given by \eqref{eq:I_f-def}.
\end{lem}
\textit{Proof of Lemma \ref{lem:one-step-backward-strong-error}}.
To establish the first inequality in \eqref{eq:one-step-backward-strong-error},
we employ
\eqref{eq:f-xy-condition},
the  H{\"o}lder inequality, Theorem
    \ref{thm:bound-backward-euler},
as well as Lemma 
\ref{lem:one-step-backward-euler-2p-bound} to deduce
\begin{equation}
 \label{eq:I_f-2p-bound}
  \begin{aligned}
   \E \big[
    |
     I _ f 
      | ^ {2p}
       \big]
  &
  =
   \E \Big[
    \Big|
     \int _ {t} ^ { t + h }
      f ( Z ( t, x; t + h ) ) -
       f (x)
       \, \dd s
        \Big| ^ {2p} 
         \Big]
%   \\
%	& 
%	 \leq
%	   h ^ { 2p - 1 }
%	    \int _ {t} ^ { t + h }
%	     \E \big[
%	      |
%	       f ( Z ( t, x; t + h ) ) -
%	        f (x)
%	         | ^ {2p}
%	          \big]
%	           \, \dd s
  \\
   &
    \leq
     c _ 1 ^ {p}
      h ^ { 2p - 1 }
       \int _ {t} ^ { t + h }
        \E \Big[
         \big(
          1 +
          |
           Z ( t, x; t + h )
            | ^
             { 2 \kappa - 2 }
   +
    | x | ^
     { 2 \kappa - 2 }
      \big) ^ {p}
       \big| 
        Z ( t, x; t + h ) - 
         x 
          \big| ^ {2p}
           \Big]
            \, \dd s
  \\
   & 
   \leq
    c _ 1 ^ {p}
     h ^ { 2p - 1 }
     \int _ {t} ^ { t + h }
      \Big(
       \E \Big[
        \big(
         1 +
          |
           Z ( t, x; t + h )
            | ^
             { 2 \kappa - 2 }
   +
    | x | ^
     { 2 \kappa - 2 }
      \big) ^ 
       {\frac{(2\kappa - 1 ) p}{\kappa -1}}
        \Big]
         \Big) ^
          { \frac{\kappa - 1}{2\kappa - 1 } }
   \\
    &
     \quad 
      \Big(
       \E \Big[
        | Z ( t, x; t + h ) - x | ^ 
         {\frac{2( 2\kappa - 1 ) p}{\kappa}}
          \Big]
           \Big) ^
            { \frac{\kappa}{2\kappa - 1 } }
             \, \dd s
%   \\
%	& 
%	 \leq
%	   K h ^ { 3p - 1 }
%	    \int _ {t} ^ { t + h }
%	     \big( 
%	      1 +
%	       | x | ^
%	        { ( 2 \kappa - 2 ) p } 
%	         \big)
%	          \big( 1 +
%	           | x | ^
%	            { 2 \kappa p }
%	             \big)
%	              \, \dd s
  \\
   & 
   \leq
    K
     \big(
       1 +
        | x | ^
         { ( 4 \kappa - 2 ) p } 
          \big)
           h ^ {3p}.
  \end{aligned}
\end{equation}
By utilizing
 \eqref{eq:f-xy-condition},
  \eqref{eq:g-xy-lipschitz}
and  \eqref{eq:X_t-moment-bound},
along with \eqref{eq:I_f-2p-bound},
we can apply  the H{\"o}lder inequality,
the moment inequality and  Lemma \ref{lem:one-step-2p-bound}
to estimate the other inequality in
   \eqref{eq:one-step-backward-strong-error}
as follows:
\begin{equation}
 \begin{aligned}
  &
  \quad
  \E \big[
   |
     X ( t, x; t + h ) -
       Z ( t, x; t + h )
      | ^ {2p}
       \big]
       \\
  &
  =
   \E \Big[
    \Big|
     \int _ {t} ^ { t + h }
      f ( X (s) ) -
       f ( Z ( t, x; t + h ) )
        \, \dd s
  +
   \int _ {t} ^ { t + h }
    g ( X (s) ) - g (x)
     \, \dd W (s)
      \Big| ^ {2p}
       \Big]
  \\
   &
   =
    \E \Big[
     \Big|
      \int _ {t} ^ { t + h }
       f ( X (s) ) - f (x)
        \, \dd s
   -
    I _ f
   +
    \int _ {t} ^ { t + h }
     g ( X (s) ) - g (x)
      \, \dd W (s)
       \Big| ^ {2p}
        \Big] 
  \\
   &
    \leq
      C _ p
       h ^ { 2p - 1 }
        \int _ {t} ^ { t + h }
         \E \big[
          |
           f ( X(s) ) - f (x)
            | ^ {2p}
             \big]
              \, \dd s
   +
    C _ p
    \E \big[
     |
      ( I _ f )
       | ^ {2p}
        \big]
   \\
	&
	\quad
	 +
	  C _ p
	   h ^ { p - 1 }
	   \int _ {t} ^ { t + h }
	    \E \big[
	     \left\|
	      g ( X(s) ) - g (x)
	       \right\| ^ {2p}
	        \big]
	         \, \dd s
   \\
    &
     \leq
       c _ 1 ^ {p}
        C _ p
         h ^ { 2p - 1 }
          \int _ {t} ^ { t + h }
           \E \Big[
            \big(
             1 +
              | X (s)
               | ^
                { 2 \kappa - 2 }
    +
     | x | ^
      { 2 \kappa - 2 }
       \big) ^ {p}
        \big| 
         X (s) -
          x 
           \big| ^ {2p}
            \Big]
             \, \dd s
   \\
	&
	\quad  
	 +
	  K 
	  \big( 
	   1 +
	   | x | ^
	    { ( 4 \kappa - 2 ) p } 
	     \big)
	      h ^ {3p}
   +
    C _ p
     h ^ { p - 1 }
     \beta _ 1 ^ {p}
      \int _ {t} ^ { t + h }
       \E \big[
        |
         X (s) - x
         | ^ {2p}
          \big] 
           \, \dd s
   \\
	& 
	\leq
	 K h ^ { 2p - 1 }
	 \int _ {t} ^ { t + h }
	  \Big(
	   \E \Big[
	    \big(
	     1 +
	      | X (s) | ^
	       { 2 \kappa - 2 }
	+
	 | x | ^
	  { 2 \kappa - 2 }
	   \big) ^
	    {\frac{(2\kappa - 1 ) p}{\kappa -1}}
	     \Big]
	      \Big) ^
	       { \frac{\kappa - 1}{2\kappa - 1 } }
   \\
	&
	 \quad 
	  \Big(
	   \E \Big[
	    | X (s) - x | ^ 
	     {\frac{2( 2\kappa - 1 ) p}{\kappa}}
	      \Big]
	       \Big) ^
	        { \frac{\kappa}{2\kappa - 1 } }
         \, \dd s
         +
	         K 
	          \big(
	           1 +
	            | x | ^
	             { ( 4 \kappa - 2 ) p } 
	              \big)
	               h ^ {3p} +
	                K 
	                 \big( 
	                  1 +
	                   | x | ^
	                    { 2 \kappa p } 
	                     \big)
	                      h ^ {2p}
   \\
	& 
	\leq
	 K
	  \big( 
	   1 +
	    | x | ^
	     { ( 4 \kappa - 2 ) p } 
	      \big)
	    h ^ {2p}.
	\end{aligned}
\end{equation}         
Thus the proof is finished.
$\square$

 \begin{thm}
 	\label{thm:global-error-backward} 
Under  Assumptions
 \ref{ass:SDEs-xy-contractive},
\ref{ass:g-xy-condition},
for any 
$ 1 
   \leq 
    p 
    \leq
      \frac{
      \lfloor p ^ {*} 
      \rfloor
      }{2\kappa - 1} $
and $ 0 < h 
 	  \leq h_1$,
the backward Euler method \eqref{eq:continuous-backward-euler}
has a strong convergence rate of order $\frac{1}{2}$ 
over infinite time,
namely,
\begin{equation}
 \begin{aligned}
    \Big(
     \E \big[
      | X _ k  -
        Z _ k | ^ {2p}
         \big]
          \Big) ^
           { \frac{1}{2p} }
   \leq
     K
     \Big(
      \E \big[
       \big( 
        1 +
         | X _ 0 | ^
          { ( 5\kappa - 3 ) p }
           \big) 
             \big]
             \Big) ^
             { \frac{1}{2p} }
              h ^ {\tfrac{1}{2}}.
 \end{aligned}
\end{equation}
\end{thm} 
\textit{Proof of Theorem \ref{thm:global-error-backward}}.
With the aid of   
Theorem \ref{thm:bound-backward-euler},
Lemmas
  \ref{lem:one-step-backward-weak-error},
  \ref{lem:one-step-backward-strong-error} 
  and according to Theorem \ref{thm:global-error},
one can straightforwardly obtain the long-time strong convergence rate of the backward Euler method \eqref{eq:continuous-backward-euler}.
$\square$

	\section{Applications:
		Strong Convergence Rate of the Projected Euler Method
		over Infinite Time}
	\label{sec:projected} 
	\noindent
 In this section, let us consider another application of the strong convergence theorem and analyze the strong convergence rate of a kind of projected Euler scheme \cite{beyn2016stochastic} under certain assumptions.
	The projected Euler method of SDEs \eqref{eq:form-SDE} proposed here is given as follows:
		\begin{equation}
			\label{eq:form-modified-Euler}
			\left\{
			\begin{aligned}
				\bar{Z} _ {k}
				& 
				:= 
				\Phi ( Z _ { k  } ),
				\\
				Z _ {k+1}
				& 
				:= 
				\bar{Z} _ k +
				h f ( \bar{Z} _ k ) +
				g( \bar{Z} _ k )
				\Delta W _ { k  },
    \quad
    Z _ 0 = X _ 0,
			\end{aligned}\right.
		\end{equation}
		where
	$ \Delta W _ { k  } := W ( t _ {k+1} ) - W ( t _ { k  } ),
	k \geq 1  $
	and the projection operator $\Phi \colon \R^d \rightarrow \R^d$ is assumed to 
    satisfy the following assumptions.
		\begin{assumption}
             \label{ass:Phi}
			Let the mapping $\Phi \colon \R^d \rightarrow \R^d$ obey 
			\begin{equation}
   \label{eq:condition-Phi}
   \begin{aligned}
   | \Phi (x) |
				 & \leq 
%				| x  |
%				 \wedge 
				 	h ^ { -\frac{1}{ 2 ( \kappa + 1 ) } } ,
      \\
				| \Phi (x) - \Phi (y) |
				 &\leq 
				 | x  - y|     
   \end{aligned}
   \end{equation}
		and $ \Phi (0)=0\in \R^d $,
		where  
 % $x \wedge y $ is denoted by the minimum of $x$ and $y$ and 
  $\kappa$ is given by \eqref{eq:f-xy-condition}.	
		\end{assumption}

Taking $y=0$ in \eqref{eq:condition-Phi} yields
\begin{equation}
| \Phi (x) |
				 \leq 
				 | x |.
\end{equation}
Below we divide the analysis of the strong convergence rate of the scheme \eqref{eq:form-modified-Euler}
 into three parts. 
 The first part is to establish the $2p$-th ($p \geq 1$) moment boundedness of the numerical solution.
	
\subsection{Bounded moments of the projected Euler method}
	\noindent   
 This part is  to establish the boundedness of the high-order moments of the projected Euler method. 
 %To achieve this goal, we will provide a specific condition regarding the diffusion coefficient of the SDEs \eqref{eq:form-SDE}, along with a necessary lemma.
	Based on Assumption \ref{ass:SDEs-xy-contractive},
 we can conclude that the diffusion coefficient of SDEs
	\eqref{eq:form-SDE} satisfies the following inequality:
	\begin{equation}
		\label{eq:g-xy-polynomial-condition}
		\left\|
		g ( x ) - g ( y ) 
		\right\| ^ { 2 } 
		\leq
		\beta_3 
		\big(
		1 +
		|  x  | ^ {  \kappa - 1 } +
		|  y | ^ {  \kappa - 1 } 
		\big) 
		|  x - y | ^ { 2 } ,
	\end{equation}
	where $\beta_3$ is a constant.
	As a consequence, we have
	\begin{equation}
		\label{eq:g-x-polynomial-condition}
		\left\|
		g ( x )  
		\right\| ^ { 2 } 
		\leq
		\beta_4 
		|  x  | ^ { \kappa + 1 } +
		c _ 6 ,
	\end{equation}
	where $ \beta _ 4 =
	\frac{ 2 \beta_3( \kappa + 3 ) }{ \kappa + 1 } $, 
	$ c _ 6 = 
	2 \left\| g(0) \right\| ^ 2 + \frac{ 2 \beta_3 (\kappa-1) }{ \kappa + 1 }$.
	In addition, 
 employing
	\eqref{eq:f-x-condition}, 
	\eqref{eq:g-x-polynomial-condition}
	and the fact that
	$ | \Phi(Z_{k-1}) | \leq  h ^ { -\frac{1}{ 2 ( \kappa + 1 ) } } $ yields
	\begin{equation}
		\label{eq:f-x-condition-projected-euler}
		\begin{aligned}
			| f ( \bar{Z} _ {k-1} ) | ^ 2 
			&
			\leq
			c _ 2
			| \bar{Z} _ {k-1} | ^ { 2 \kappa} +
			c _ 3 
			\leq 
			c _ 2 h ^ {-1} +
			c _ 3
		\end{aligned}
	\end{equation}
	and
	\begin{equation}
		\label{eq:g-x-condition-projected-euler}
		\begin{aligned}
			\left\| 
			g ( \bar{Z} _ {k-1} ) 
			\right\| ^ 2 
			&
			\leq
			\beta_4 
			| \bar{Z} _ {k-1} | ^ { \kappa + 1 } +
			c _ 6 
			\leq 
			\beta_4 h ^ {- \frac{1}{2}} +
			c _ 6.
		\end{aligned}
	\end{equation}

	Next,
	we shall  establish the moment bounds for the projected Euler scheme.

	\begin{thm}
		\label{thm:bound-projected-Euler}
		Let Assumptions
		\ref{ass:SDEs-xy-contractive},\ref{ass:Phi}
		hold.
		For any
		$ 1 
		\leq 
		p 
		\leq
  \lfloor
		p ^ {*}
  \rfloor$ 
		and 
		$ 0 < h 
		\leq h_2 
		:= \min
		\{
		h _ 1,
		\frac{1}{ 2p ( \alpha _ 1 - \epsilon_1)}
		\} $,
		one has
		\begin{equation}
			\label{eq:bound-projected-Euler}
			\begin{aligned}
				\E \big[
				| Z _ k | ^ {2p}
				\big]
				&
				\leq 
				K
				\E \big[
				\big( 
				1 +
				| X _ 0 | ^ {2p}
				\big)
				\big].
			\end{aligned}
		\end{equation}
	\end{thm}
	\textit{Proof of Theorem \ref{thm:bound-projected-Euler}}.
	From 
	\eqref{eq:form-modified-Euler},
	one can derive 
	\begin{equation}
		\begin{aligned}
			1 +
			| Z _ k | ^ 2
			&
			=
			1 +
			| \bar{Z} _ {k-1} +
			f( \bar{Z} _ {k-1} ) h +
			g( \bar{Z} _ {k-1})
			\Delta W _ { k  - 1 } | ^ 2
			\\
			&
			=
			1 +
			| \bar{Z} _ {k-1} | ^ 2 +
			2\left\langle
			\bar{Z} _ {k-1}, f( \bar{Z} _ {k-1} ) h
			\right \rangle +
			| g( \bar{Z} _ {k-1} )
			\Delta W _ { k  - 1 } | ^ 2 
			\\
			&
			\quad
			+
			2\left\langle
			\bar{Z} _ {k-1} +
			f( \bar{Z} _ {k-1} ) h,
			g( \bar{Z} _ {k-1} )
			\Delta W _ { k  - 1 }
			\right \rangle +
			|f( \bar{Z} _ {k-1} ) h | ^ 2.
		\end{aligned}
	\end{equation}
	Combining this with \eqref{eq:f-x-condition-projected-euler} shows
	\begin{equation}
		\label{eq:p-power-inequality-projected}
		\E [
		\big(
		1 +
		| Z _ k | ^ 2
		\big) ^ p |
		\mathcal{F} _ { t _  {k-1} } ] 	
		\leq
		\E [
		( 1 + \xi ) ^ p |
		\mathcal{F} _ { t _  {k-1} } ] 
		\big( 
		1 + 
		| \bar{Z} _ {k-1} | ^ 2 
		\big) ^ p +
		Kh
		\big( 
		1 + 
		| \bar{Z} _ {k-1} | ^ 2 
		\big) ^ {p-1},
	\end{equation}
	where we denote
	\begin{center}
		$ \xi
		:= \tfrac{
			2\left\langle
			\bar{Z} _ {k-1}, f( \bar{Z} _ {k-1} ) h
			\right \rangle +
			|g( \bar{Z} _ {k-1} )
			\Delta W _ { k  - 1 } | ^ 2  +
			2\left\langle
			\bar{Z} _ {k-1} +
			f( \bar{Z} _ {k-1} ) h,
			g( \bar{Z} _ {k-1} )
			\Delta W _ { k  - 1 }
			\right \rangle}
		{ 1 +
			|\bar{Z} _ {k-1} | ^ 2 }
		:=
		\xi _ 1 +
		\xi _ 2 +
		\xi _ 3 $.  
	\end{center}
%	If we obtain 
%	\begin{center}
%		$ \E [ ( 1 + \xi ) ^ p | \mathcal{F} _ { t _  {k-1} } ] = 1 - \gamma h$,  
%	\end{center}
%	where $ \gamma > 0 $ and independent of $ h $.
%	With the help of Lemma \ref{lem:one-step-2p-relation}, we can prove Theorem \ref{thm:bound-projected-Euler}.
%	Next, we will verify these ideas.
In what follows four cases are considered.

	\textit{Case 1}.
	When $ p = 1 $,
	one obtains that 
	\begin{equation}
		\label{eq:1-power-inequality-projected-relation}
		\begin{aligned}
			\E \big[
			\big( 
			1 +
			| Z _ k | ^ 2 
			\big) |
			\mathcal{ F } _ { t _ {k-1} } 
			\big]
			&
			\leq
			\E [
			( 1 + \xi ) |
			\mathcal{F} _ { t _ {k-1} } ]
			( 1 +
			|\bar{Z} _ {k-1}| ^ 2 ) +
			Kh
			\\
			&
			= 
			\E [
			( 1 +
			\xi _ 1 +
			\xi _ 2 +
			\xi _ 3 ) |
			\mathcal{F} _ { t _ {k-1} } ]
			\big( 
			1 +
			|\bar{Z} _ {k-1}| ^ 2 
			\big) +
			Kh.
		\end{aligned}
	\end{equation}  
	Utilizing the properties of Brownian motion,
	we have the following:
	\begin{equation}
		\begin{aligned}
			\E [
			\xi _ 1 +
			\xi _ 2  |
			\mathcal{F} _ { t _  {k-1} } ]
			&
			=
			\E \bigg[
			\frac{ 
				2
				\left\langle
				\bar{Z} _ {k-1}, f( \bar{Z} _ {k-1} ) h
				\right \rangle +
				|g( \bar{Z} _ {k-1} )
				\Delta W _ { k  - 1 } | ^ 2 }
			{ 1 +
				| \bar{Z} _ {k-1} | ^ 2 }
			\bigg|
			\mathcal{F} _ { t _ {k-1} }
			\bigg]
			\\
			&
			=
			\frac{ 
				2
				\left\langle
				 \bar{Z} _ {k-1}, f( \bar{Z} _ {k-1} ) h
				\right \rangle +
				\left\|g( \bar{Z} _ {k-1} )
				\right\| ^ 2 h}
			{ 1 +
				| \bar{Z} _ {k-1} | ^ 2 }  
		\end{aligned}
	\end{equation} 
	and
	\begin{equation}
		\begin{aligned}
			\E [
			\xi _ 3 |
			\mathcal{F} _ { t _ {k-1} } ]
			&
			=
			\E \bigg[
			\frac{
				2
				\left\langle
				 \bar{Z} _ {k-1} +
				f( \bar{Z} _ {k-1} ) h,
				g( \bar{Z} _ {k-1} )
				\Delta W _ { k  - 1 }
				\right \rangle}
			{ 1 +
				| \bar{Z} _ {k-1} | ^ 2 }
			\bigg|
			\mathcal{F} _ { t _ {k-1} }
			\bigg] 
			&
			=
			0.
		\end{aligned}
	\end{equation}  
	Taking these two estimates into account,
	we can derive from \eqref{eq:1-power-inequality-projected-relation}
	that 
	\begin{equation}
		\label{eq:1-power-coefficient-bound}
		\begin{aligned}
			\E [ 
			\xi |
			\mathcal{F} _ { t _ {k-1}} ]
			=
			\frac{ 
				2
				\left\langle
				 \bar{Z} _ {k-1}, f( \bar{Z} _ {k-1} ) h
				\right \rangle +
				\left\|g( \bar{Z} _ {k-1} )
				\right\| ^ 2 h}
			{ 1 +
				| \bar{Z} _ {k-1} | ^ 2 }.  
		\end{aligned}
	\end{equation}
	By  \eqref{eq:SDEs-x-contractive}
	and     \eqref{eq:condition-Phi},
	we acquire
	\begin{equation}
		\begin{aligned}
			\E \big[
			\big( 
			1 +
			| Z _ k | ^ 2 
			\big) 
			\big]
			&
			\leq
			( 1 -
			2 
			( \alpha_1 -
			\epsilon_1 ) 
			h )
			\E \big[ 
			\big( 
			1 +
			| Z _ { k - 1 } |^2
			\big)
			\big] +
			Kh .
		\end{aligned}
	\end{equation}
	Then from Lemma \ref{lem:one-step-2p-relation} 
	we can see that the scheme \eqref{eq:bound-projected-Euler}
	is evident when $ p = 1 $.  
	
	\textit{Case 2}. When $ p = 2 $, 
	one sees that
	\begin{equation}
		\begin{aligned}
			\E \Big[
			\big(
			1 +
			| Z _ k | ^ 2
			\big) ^ 2 
			\big|
			\mathcal{F} _ { t _ {k-1} }
			\Big]
			&
			\leq
			\E \big[
			( 1 + \xi ) ^ 2 |
			\mathcal{F} _ { t _ {k-1} }
			\big]
			\big( 
			1 + 
			| \bar{Z} _ {k-1} | ^ 2
			\big) ^ 2 +
			K h 
			\big( 
			1 +
			| \bar{Z} _ {k-1} | ^ 2 
			\big) 
			\\
			&
			=
			\E \big[
			( 1 +
			2 \xi +
			\xi ^ 2 ) |
			\mathcal{F} _ { t _ {k-1} }
			\big] 
			\big( 
			1 +
			| \bar{Z} _ {k-1} | ^ 2 
			\big) ^ 2 +
			K h
			\big( 
			1 +
			| \bar{Z} _ {k-1} | ^ 2
			\big)
		\end{aligned}
	\end{equation}   
	and 
	\begin{equation}
		\label{eq:2-power-projected-coefficient-decomposition}
		\begin{aligned}
			\E \big[
			\xi ^ 2  |
			\mathcal{F} _ { t _ {k-1} }
			\big]
			&
			=
			\E \big[
			( \xi _ 1 +
			\xi _ 2 +
			\xi _ 3 ) ^ 2 
			\big|
			\mathcal{F} _ { t _ {k-1} }
			\big]
			\\
			&
			=
			\E \big[
			\big(
			\xi _ 1 ^ 2 +
			\xi _ 2 ^ 2 +
			\xi _ 3 ^ 2 +
			2 \xi _ 1 \xi _ 2 +
			2 \xi _ 1 \xi _ 3 +
			2 \xi _ 2 \xi _ 3 
			\big) 
			\big|
			\mathcal{F} _ { t _ {k-1} }
			\big].
		\end{aligned}
	\end{equation}
	Next, 
 we will provide estimates for each term in
 \eqref{eq:2-power-projected-coefficient-decomposition}.
	By  \eqref{eq:f-x-condition-projected-euler},
	the first item on the right-hand side of \eqref{eq:2-power-projected-coefficient-decomposition}
	is bounded as:
	\begin{equation}
		\label{eq:first-item-coefficient-estimate}
		\begin{aligned}
			\E \big[
			\xi _ 1 ^ 2  |
			\mathcal{F} _ { t _ {k-1} }
			\big]
			&
			\leq
			\E \bigg[
			\tfrac{
				4| \bar{Z} _ {k-1} | ^ 2
				| f ( \bar{Z} _ {k-1}) h | ^ 2 }
			{ \big(
				1 +
				| \bar{Z} _ {k-1} | ^ 2 
				\big) ^ 2 }
			\bigg |
			\mathcal{F} _ {t _ {k-1} }
			\bigg]
%			\\
%			&
			\leq
			\frac{
				4 | \bar{Z} _ {k-1} | ^ 2
				\big(
				c _ 2 h + 
				c _ 3 h ^ 2 
				\big) }
			{ \big( 
				1 +
				| \bar{Z} _ {k-1} | ^ 2 
				\big) ^ 2 }
%			\\
%			&
			\leq
			\frac{ K h }
			{ 1 +
				| \bar{Z} _ {k-1} | ^ 2 }.
		\end{aligned}
	\end{equation} 
	For the second item on the right-hand side of \eqref{eq:2-power-projected-coefficient-decomposition},
	using \eqref{eq:g-x-condition-projected-euler} gives
	\begin{equation}
		\label{eq:second-item-coefficient-estimate}
		\begin{aligned}
			\E \big[
			\xi _ 2 ^ 2  |
			\mathcal{F} _ { t _ {k-1} }
			\big]
%			&
%			=
			\E \bigg[
			\frac{
				| g ( \bar{Z} _ {k-1})
				\Delta W _ { k  - 1 } | ^ 4 }
			{ \big (
				1 +
				| \bar{Z} _ {k-1} | ^ 2 
				\big ) ^ 2 }
			\bigg |
			\mathcal{F} _ { t _ {k-1} }
			\bigg]
%			\\
%			&
			\leq  
			\frac{
				3 
				\big(
				\beta_4 h ^ {- \frac{1}{2}} +
				c _ 6
				\big) ^ 2
				h ^ 2 }
			{ \big(
				1 +
				| \bar{Z} _ {k-1} | ^ 2 
				\big ) ^ 2 }
%			\\
%			&
			\leq
			\frac{ K h }{	1 +
				| \bar{Z} _ {k-1} | ^ 2 } .
		\end{aligned}
	\end{equation}
	Employing \eqref{eq:SDEs-x-contractive}
	and \eqref{eq:g-x-lipschitz}
	to bound the third and  fourth items on the right-hand side of  \eqref{eq:2-power-projected-coefficient-decomposition} as follows:
	\begin{equation}
		\label{eq:third-item-coefficient-estimate}
		\begin{aligned}
			\E \big[
			\xi _ 3 ^ 2  |
			\mathcal{F} _ { t _ {k-1} }
			\big]
			&
		\leq
			\E \bigg[
			\frac{
				4 |
				 \bar{Z} _ {k-1} +
				f ( \bar{Z} _ {k-1}) h | ^ 2
				| g ( \bar{Z} _ {k-1} )
				\Delta W _ { k  - 1 } | ^ 2 }
			{ \big(
				1 +
				| \bar{Z} _ {k-1} | ^ 2 
				\big) ^ 2 }
			\bigg |
			\mathcal{F} _ { t _ {k-1} }
			\bigg]
			\\
			& 
			\leq 
			\frac{ 4h
				\left\|g ( \bar{Z} _ {k-1} )
				\right\| ^ 2  }{	1 +
				| \bar{Z} _ {k-1} | ^ 2 } +
			\frac{ K h }{	1 +
				| \bar{Z} _ {k-1} | ^ 2 } 
		\end{aligned}
	\end{equation}
	and
	\begin{equation}
		\label{eq:fourth-item-coefficient-estimate}
		\begin{aligned}
			\E \big[
			2 \xi _ 1
			\xi _ 2 |
			\mathcal{F} _ { t _ {k-1} }
			\big]
			&
			=
			\E \bigg[
			\frac{
				4 \left\langle
				 \bar{Z} _ {k-1}, 
				f ( \bar{Z} _ {k-1}) h
				\right \rangle
				| g ( \bar{Z} _ {k-1} )
				\Delta W _ { k  - 1 } 
				| ^ 2 }
			{ \big( 
				1 +
				| \bar{Z} _ {k-1} | ^ 2 
				\big) ^ 2 }
			\bigg |
			\mathcal{F} _ { t _ {k-1} }
			\bigg] 
%			\\
%			& 
			\leq 
			\frac{ K h }
			{ 1 +
				| \bar{Z} _ {k-1} | ^ 2 }.
		\end{aligned}
	\end{equation}
Due to the properties of Brownian motion, 
	we know
	\begin{center}
		$ \E \big[
		2 \xi _ 1 \xi _ 3 |
		\mathcal{F} _ { t _ {k-1} }
		\big]
		=
		\E \big[
		2 \xi _ 2 \xi _ 3 |
		\mathcal{F} _ { t _ {k-1} }
		\big]
		=
		0 $.  
	\end{center}
	Combining \eqref{eq:first-item-coefficient-estimate},
	 \eqref{eq:second-item-coefficient-estimate},
	\eqref{eq:third-item-coefficient-estimate} 
	with
	\eqref{eq:fourth-item-coefficient-estimate}, we obtain
	\begin{equation}
		\begin{aligned}
			\label{eq:2-power-coefficient-bound}
			\E \big[
			\xi  ^ 2
			\big|
			\mathcal{F} _ { t _ {k-1} }
			\big]
			\leq
			\frac{ 4h
				\left\|g ( \bar{Z} _ {k-1} )
				\right\| ^ 2  }{	1 +
				| \bar{Z} _ {k-1} | ^ 2 } +
			\frac{ K h }
			{ 1 +
				| \bar{Z} _ {k-1} | ^ 2 }.    
		\end{aligned}
	\end{equation}
	Next, 
	by \eqref{eq:SDEs-x-contractive} 
	and \eqref{eq:condition-Phi},
	we have
	\begin{equation}
		\begin{aligned}
			\E \Big[
			\big
			( 1 +
			| Z _ k | ^ 2
			\big ) ^ 2 
			\big|
			\mathcal{F} _ { t _ {k-1} }
			\Big]
			&
			\leq
			\Big(
			1 +
			\tfrac{4 \left\langle
				 \bar{Z} _ {k-1}, 
				f ( \bar{Z} _ {k-1}) h
				\right \rangle +
				6h
				\left\|
				g ( \bar{Z} _ {k-1} )
				\right\| ^ 2 }{	1 +
				| \bar{Z} _ {k-1} | ^ 2 }
			\Big)
			\big(
			1 +
			| \bar{Z} _ {k-1} | ^ 2
			\big) ^ 2 +
			K h
			\big( 
			1 +
			| \bar{Z} _ {k-1} | ^ 2 
			\big)
			\\
			&
			\leq
			( 1 -
			4 ( \alpha_1 -
			\epsilon_1 ) 
			h)
			\big( 
			1 +
			| Z _ { k - 1 } | ^ 2
			\big ) ^ 2  +
			K h
			\big( 
			1 +
			| Z _ { k - 1 } | ^ 2
			\big).  	               
		\end{aligned}
	\end{equation}
	Thanks to Lemma \ref{lem:one-step-2p-relation}, 
	the inequality \eqref{eq:bound-projected-Euler}
	holds when $ p = 2 $. 
	
	\textit{Case 3}. For $ p = 3 $,
	we make a further decomposition as follows:
	\begin{equation}
		\begin{aligned}
			\E \Big[
			\big
			( 1 +
			| Z _ k | ^ 2
			\big ) ^ 3 
			\big|
			\mathcal{F} _ { t _ {k-1} }
			\Big] 
			&
			\leq
			\E \big[
			( 1 +
			\xi ) ^ 3 
			\big|
			\mathcal{F} _ { t _ {k-1} }
			\big] 
			\big
			( 1 +
			| \bar{Z} _ {k-1} | ^ 2
			\big) ^ 3 +
			K h 
			\big(
			1 +
			| \bar{Z} _ {k-1} | ^ 2 
			\big) ^ 2
			\\
			&
			=  
			\E \big[
			( 1 +
			3 \xi +
			3 \xi ^ 2 +
			\xi ^ 3 ) 
			\big|
			\mathcal{F} _ { t _ {k-1} }
			\big] 
			\big
			( 1 +
			| \bar{Z} _ {k-1} | ^ 2
			\big ) ^ 3 +
			K h
			\big(
			1 +
			| \bar{Z} _ {k-1} | ^ 2 
			\big) ^ 2
		\end{aligned}
	\end{equation} 
	and
	\begin{equation}
		\label{eq:3-power-projected-coefficient-decomposition}
		\begin{aligned}
			\E \big[
			\xi ^ 3 
			\big|
			\mathcal{F} _ { t _ {k-1} }
			\big]
			&
			=
			\E \big[
			( \xi _ 1 +
			\xi _ 2 +
			\xi _ 3 ) ^ 3 
			\big|
			\mathcal{F} _ { t _ {k-1} }
			\big] 
			\\
			&
			=
			\E \big[
			\big(
			\xi _ 2 ^ 3 +
			3 \xi _ 1 ^ 2 \xi _ 2 +
			3 \xi _ 2 \xi _ 3 ^ 2 +
			\xi _ 1 ^ 3 +
			3 \xi _ 1 \xi _ 2 ^ 2 +
			3 \xi _ 1 \xi _ 3 ^ 2 +
			3( \xi _ 1 +
			\xi _ 2 ) ^ 2 \xi _ 3 +
			\xi _ 3 ^ 3 
			\big) 
			\big|
			\mathcal{F} _ { t _ {k-1} }
			\big].
		\end{aligned}
	\end{equation}
	Next,
	we will bound these terms separately. 
	By applying   
	\eqref{eq:g-x-condition-projected-euler}, 
	we obtain
	\begin{equation}
		\label{eq:3-power-first-item}	
		\begin{aligned}
			\E \big[
			\xi _ 2 ^ 3 
			\big|
			\mathcal{F} _ { t _ {k-1} }
			\big]
			&
			=
			\E \bigg[
			\frac{
				| g ( \bar{Z} _ {k-1} )
				\Delta W _ { k  - 1 } | ^ 6  }
			{ \big( 
				1 +
				| \bar{Z} _ {k-1} | ^ 2
				\big) ^ 3 }
			\bigg |
			\mathcal{F} _ { t _ {k-1} }
			\bigg]
			\leq  
			\frac{
				15 
				\big(
				\beta_4
				h ^{ -\frac{1}{2} } +
				c _ 6 
				\big) ^ 3
				h ^ 3 }
			{ \big( 
				1 +
				| \bar{Z} _ {k-1} | ^ 2
				\big) ^ 3 }
			\leq 
			\frac{ K h}
			{ 1 +
				| \bar{Z} _ {k-1} | ^ 2 }.
		\end{aligned}
	\end{equation} 
	With the help of \eqref{eq:f-x-condition-projected-euler}
	and \eqref{eq:g-x-condition-projected-euler},
	the second and third items on the right-hand side of  \eqref{eq:3-power-projected-coefficient-decomposition}
	can be estimated as
	\begin{equation}
		\begin{aligned}
			\E \big[
			3 \xi _ 1 ^ 2 \xi _ 2 
			\big|
			\mathcal{F} _ { t _ {k-1} }
			\big]
			&
			=
			\E \bigg[
			\frac{
				12 | \bar{Z} _ {k-1} | ^ 2
				| f ( \bar{Z} _ {k-1} ) h | ^ 2
				| g ( \bar{Z} _ {k-1} )
				\Delta W _ { k  - 1 } | ^ 2}
			{ \big
				( 1 +
				| \bar{Z} _ {k-1} | ^ 2
				\big ) ^ 3 }
			\bigg |
			\mathcal{F} _ { t _ {k-1} }
			\bigg]
			\\
			&
			\leq
			\frac{
				12 
				| \bar{Z} _ {k-1} | ^ 2
				\big(
				c _ 2
				h ^ {-1} +
				c _ 3
				\big)
				\big(
				\beta_4
				h ^{ -\frac{1}{2} } +
				c _ 6 
				\big)
				h ^ 3}
			{ \big(
				1 +
				| \bar{Z} _ {k-1} | ^ 2
				\big ) ^ 3 }      
			\\
			&
			\leq
			\frac{ Kh }
			{ 1 +
				| \bar{Z} _ {k-1} | ^ 2 } 
		\end{aligned}
	\end{equation}
	and
	\begin{equation}
		\begin{aligned}
			\E \big[
			3 \xi _ 2 \xi _ 3 ^ 2 
			\big|
			\mathcal{F} _ { t _ {k-1} }
			\big]
			&
			\leq
			\E \bigg[
			\frac{
				12
				| \bar{Z} _ {k-1} +
				f( \bar{Z} _ {k-1} ) h | ^ 2
				| g ( \bar{Z} _ {k-1} )
				\Delta W _ { k  - 1 } | ^ 4 }
			{ \big
				( 1 +
				| \bar{Z} _ {k-1} | ^ 2
				\big) ^ 3 }
			\bigg |
			\mathcal{F} _ { t _ {k-1} }
			\bigg]
			\leq      
			\frac{ Kh }
			{ 1 +
				| \bar{Z} _ {k-1} | ^ 2 } .
		\end{aligned}
	\end{equation}
In view of \eqref{eq:SDEs-x-contractive},
	it is not difficult to observe that
	\begin{equation}
		\label{eq:3-power-second-item}
		\begin{aligned}
			\E \big[
			\big(
			\xi _ 1 ^ 3 +
			3 \xi _ 1 \xi _ 2 ^ 2 +
			3 \xi _ 1 \xi _ 3 ^ 2 
			\big) 
			\big|
			\mathcal{F} _ { t _ {k-1} }
			\big]
			&
			\leq
			\frac{ Kh }
			{ 1 +
				| \bar{Z} _ {k-1} | ^ 2 } .		
		\end{aligned}
	\end{equation}	
By applying properties of the Brownian motion,
	we know
	\begin{center}
		$ \E \big[
		\xi _ 3 ^ 3 
		\big|
		\mathcal{F} _ { t _ {k-1}}
		\big]
		=
		\E \big[
		3
		( \xi _ 1 +
		\xi _ 2 ) ^ 2
		\xi _ 3|
		\mathcal{F} _ { t _ {k-1} }
		\big]
		=
		0 $.
	\end{center}
	Then, using \eqref{eq:3-power-first-item}-\eqref{eq:3-power-second-item}, 
	we are able to obtain
	\begin{equation}
		\label{eq:3-power-coefficient-bound}
		\E \big[
		\xi  ^ 3 
		\big|
		\mathcal{F} _ { t _ {k-1} }
		\big]
		\leq
		\tfrac{ Kh }
		{ 1 +
			| \bar{Z} _ {k-1} | ^ 2 } .	
	\end{equation}        
	By  \eqref{eq:SDEs-x-contractive} 
	and  \eqref{eq:condition-Phi}, 
	we see
	\begin{equation}
		\begin{aligned}
			\E \Big[
			\big
			( 1 +
			| Z _ k | ^ 2
			\big ) ^ 3 
			\big|
			\mathcal{F} _ { t _ {k-1} }
			\Big]
			&
			\leq
			\Big(
			1 +
			\tfrac{6 
				\left\langle
				 \bar{Z} _ {k-1}, 
				f ( \bar{Z} _ {k-1}) h
				\right \rangle +
				15h
				\left\|
				g ( \bar{Z} _ {k-1} )
				\right\| ^ 2 }{	1 +
				| \bar{Z} _ {k-1} | ^ 2 }
			\Big)
			\big(
			1 +
			| \bar{Z} _ {k-1} | ^ 2
			\big) ^ 3 
   \\
   & \quad
   +
			K h
			\big(
			1 +
			| \bar{Z} _ {k-1} | ^ 2 
			\big) ^ 2 
			\\
			&
			\leq
			( 1 -
			6 ( \alpha_1 -
			\epsilon_1 ) 
			h )
			\big( 
			1 +
			| Z _ { k - 1 }  | ^ 2
			\big ) ^ 3 + 
			K h
			\big(
			1 +
			| Z _ { k - 1 }  | ^ 2 
			\big) ^ 2.                  
		\end{aligned}
	\end{equation}	
	By applying Lemma \ref{lem:one-step-2p-relation}, 
	the proof of the inequality \eqref{eq:bound-projected-Euler}
	is thus completed  for the case when $ p = 3 $.
	
	%For $ p \geq 4 $, 
	%according to the inequality \eqref{eq:p-power-inequality-projected}, 
	%we have
	%\begin{equation}
	%	\E \big[
	%	 \xi _ 1 ^ 2  
	%	  \big|
	%	   \mathcal{F} _ { t _  {k-1} } 
	%	    \big] 
	%   \leq
	%	\frac{ K h}
	%	 { 1 +
		%	  | \bar{Z} _ {k-1} | ^ 2 },
	%\end{equation} 
	%\begin{equation}
	%	\E \big[
	%	 \xi _ 2 ^ 2  
	%	  \big|
	%	   \mathcal{F} _ { t _  {k-1} } 
	%	    \big] 
	%   \leq
	%	\frac{ K h}
	%	 { 1 +
		%	  | \bar{Z} _ {k-1} | ^ 2 },
	%\end{equation} 
	%\begin{equation}
	%	\E \big[
	%	 \xi _ 3 ^ 4  |
	%	  \mathcal{F} _ { t _  {k-1} } 
	%	   \big] 
	%   \leq
	%	\frac{ K h}
	%	 { 1 +
		%	  | \bar{Z} _ {k-1} | ^ 2 }.
	%\end{equation} 
	%Thus the inequality 
	%\eqref{eq:p-power-inequality-projected} implies 
\textit{Case 4}.	
 For $ p \geq 4 $, 
	based on the previous analysis, 
	we have
	\begin{equation}
		\label{eq:p-power-inequality-projected-relation}
		\begin{aligned}
			\E \big[
			\big(
			1 +
			| Z _ k | ^ 2
			\big) ^ p
			\big|
			\mathcal{F} _ { t _  {k-1} } 
			\big] 	
			&
			\leq
			\E \Big[
			\big(
			1 + 
			p \xi +
			\tfrac{ p (p-1) }{2}
			\xi ^ 2 +
			\tfrac{ p (p-1) (p-2) }{6}
			\xi ^ 3 \big) \big|
			\mathcal{F} _ { t _ {k-1}} 
			\Big]
			\\
			&
			\quad  
			\big( 
			1 + 
			| \bar{Z} _ {k-1} | ^ 2 
			\big) ^ p 		
			+
			Kh
			\big( 
			1 + 
			| \bar{Z} _ {k-1} | ^ 2 
			\big) ^ {p-1}.
		\end{aligned}
	\end{equation}
	Taking \eqref{eq:1-power-coefficient-bound}, 
	\eqref{eq:2-power-coefficient-bound} 
	and \eqref{eq:3-power-coefficient-bound} 
	into account and using \eqref{eq:SDEs-x-contractive},
	 \eqref{eq:condition-Phi}
	we derive from 
	\eqref{eq:p-power-inequality-projected-relation} 
	that
	\begin{equation}
		\begin{aligned}
			\E \big[
			\big(
			1 +
			| Z _ k | ^ 2
			\big) ^ p
			\big|
			\mathcal{F} _ { t _  {k-1} } 
			\big] 	
			&
			\leq
			\bigg(
			1 +
			\tfrac{ 2p 
				\left\langle
				 \bar{Z} _ {k-1}, 
				f ( \bar{Z} _ {k-1}) h
				\right \rangle +
				p ( 2p - 1 )
				\left\|
				g ( \bar{Z} _ {k-1} )
				\right\| ^ 2 }{	1 +
				| \bar{Z} _ {k-1} | ^ 2 }
			\bigg)	           
			\\
			&
			\quad
		\times	
   \big(
			1 +
			| \bar{Z} _ {k-1} | ^ 2
			\big) ^ p
			+
			K h
			\big(
			1 +
			| \bar{Z} _ {k-1} | ^ 2 
			\big) ^ {p-1} 
			\\
			&
			\leq
			( 1 -
			2p 
			( \alpha_1 -
			\epsilon_1 ) 
			h )
			\big(
			1 +
			| Z _ { k - 1 } | ^ 2
			\big) ^ p +
			K h
			\big(
			1 +
			| Z _ { k - 1 }| ^ 2 
			\big) ^ {p-1}.
		\end{aligned}
	\end{equation}
	Therefore,  \eqref{eq:bound-projected-Euler} 
	is evident for integer $p$ by applying Lemma 
	\ref{lem:one-step-2p-relation}. 
	For non-integer $p$, we can bound these  terms using  the Young inequality, 
	which completes the proof.   
	$\square$
	
	\subsection{Strong convergence rate of the projected Euler method}
	\noindent  
	In light of the fundamental strong convergence theorem,
	Theorem 
	\ref{thm:global-error},
	we need to verify the local weak and strong errors in \eqref{eq:one-step-weak-error}
	and 
	\eqref{eq:one-step-strong-error}.
	For the purpose of the strong convergence analysis,
	we will start from the following auxiliary lemmas. 
	
	\begin{lem}
		\label{lem:one-step-2p-projected-bound}
		Suppose Assumption 
		\ref{ass:SDEs-xy-contractive} 
		hold.
		For 
		$ 1 
		\leq 
		p 
		<
		\frac{p ^ {*} }{\kappa} $,
		one obtains that
		\begin{equation}
			\label{eq:one-step-2p-projected-bound}
			\begin{aligned}
				\E \big[
				| X ( h ) - X _ 0 | ^ {2p}
				\big]
				\leq
				C _ 8
				\E \big[
				\big(
				1 +
				| X _ 0 | ^ { 2 \kappa p } 
				\big)
				\big]
				h ^ {p}.
			\end{aligned}
		\end{equation}
		where $ C _ 8 $  do not depend on $h,t$.	
	\end{lem}
	\textit{Proof of Lemma 
		\ref{lem:one-step-2p-projected-bound}.}
	Based on
	\eqref{eq:f-x-condition}
	and \eqref{eq:g-x-condition-projected-euler},
	we can use the H{\"o}lder inequality, the moment inequality and \eqref{eq:X_t-moment-bound}
	to obtain
	\begin{equation}
		\begin{aligned}
			\E \big[
			| X ( h ) - X _ 0 | ^ {2p}
			\big]
			&
			=
			\E \Big[
			\Big|
			\int _ {0} ^ {h}
			f ( X (s) )
			\, \dd s +
			\int _ {0} ^ {h}
			g ( X (s) )
			\, \dd  W(s)
			\Big| ^ {2p}
			\Big]
			%   \\
			%	&
			%	\leq 
			%	 C _ p 
			%	  h ^ { 2p - 1 }
			%	   \int _ {0} ^ {h}
			%		\E \big[
			%	 	 | f ( X (s) ) | ^ {2p}
			%		  \big]
			%		   \, \dd s +
			%		    C _ p
			%		     h ^ { p - 1 }
			%		      \int _ {0} ^ {h}
			%		       \E \big[
			%		        \left\|
			%		         g ( X (s) )
			%		          \right\|^ {2p}
			%		           \big]
			%		            \, \dd  s
			\\
			&
			\leq 
			C _ p
			h ^ { 2p - 1 }
			\int _ {0} ^ {h}
			\E \Big[
			\big(
			c _ 3 +
			c _ 2
			| X (s) | ^ { 2 \kappa }
			\big) ^ {p}
			\Big]
			\, \dd s +
			C _ p
			h ^ { p - 1 }
			\int _ {0} ^ {h}
			\E \Big[
			\big(
			c _ 6 +
			\beta _ 4
			| X (s) | ^ { \kappa + 1 }
			\big) ^ {p}
			\Big]
			\, \dd s
			\\
			& 
			\leq
			C _ 8
			\E \big[
			\big(
			1 +
			| X _ 0 | ^ { 2 \kappa p } 
			\big)
			\big]
			h ^ {p}.
		\end{aligned}
	\end{equation}
	Thus the proof is finished.
	$\square$

	\begin{lem}
		\label{lem:error-of-x-projected}
		For $ h \in ( 0, 1 ] $ and $ q \geq 1 $,
		consider the mapping
		$ \R ^ d \ni x \mapsto \Phi(x) \in \R ^ d $
		which satisfies
		\begin{equation}
			\label{eq:error-of-x-projected}
			| x - \Phi(x) |
			\leq
		2 
			\big( 1 +
			| x | ^ { q + 1 } 
			\big)
			h ^ { \frac{ q  }{  2 (  \kappa + 1 ) } }.
		\end{equation} 
	\end{lem}
	\textit{Proof of Lemma \ref{lem:error-of-x-projected} }.
	We will divide the analysis into two cases: 
	$ | x | \leq h ^ { -\frac{1}{  2 (  \kappa + 1 ) } } $
	and $ | x | > h ^ { -\frac{1}{  2 (  \kappa + 1 ) } } $
	.
	
	For $ | x | \leq h ^ { -\frac{1}{  2 (  \kappa + 1 ) } } $,
	we have $ x = \Phi(x) $, 
	thus \eqref{eq:error-of-x-projected} is evident.
	
	For $ | x | > h ^ { -\frac{1}{  2 (  \kappa + 1 ) } } $,
	by using the elementary inequality, we can easily obtain
	\begin{equation}
		\begin{aligned}
			| x - \Phi(x) |  
			=
		\Big| x-
		h ^ { -\frac{1}
			{ 2 ( \kappa + 1 ) } } 
		\frac{x}{|x|}		
		\Big|  
			\leq
			| x | +
			\big| 
			h ^
			{-\frac{1} { 2 (  \kappa + 1 ) }} 
			\big| 
			\leq
		2 
			\big( 
			1 +
			| x | ^ { q + 1 } 
			\big)
			h ^ { \frac{ q  }{ 2 (  \kappa + 1 ) } }.
		\end{aligned}
	\end{equation}
	Thus the proof is accomplished.
	$\square$
	\begin{lem}
		\label{lem:one-step-projected-weak-error}
		Under the Assumptions \ref{ass:SDEs-xy-contractive} and \ref{ass:Phi}, 
		we obtain 
		\begin{equation}
			\label{eq:one-step-projected-weak-error}
			\begin{aligned}
				\big|
				\E \big[
				Z _ E ( t, x; t + h ) -
				Z( t, x; t + h ) 
				\big]
				\big|
				&
				\leq
				K 
				\big( 
				1 +
				| x | ^
				{  3 \kappa + 4 } 
				\big)
				h ^ { \frac{3}{2} },
				\\
				\big|
				\E \big[
				X ( t, x; t + h ) -
				Z( t, x; t + h )
				\big] 
				\big|
				&
				\leq
				K 
				\big(
				1 +
				| x | ^
				{ 3 \kappa + 4 } 
				\big)
				h ^ { \frac{3}{2} }.
			\end{aligned}
		\end{equation}   
	\end{lem} 
	\textit{Proof of Lemma \ref{lem:one-step-projected-weak-error} }.
	We can employ 
	\eqref{eq:f-xy-condition}, 
	the  H{\"o}lder inequality  and Lemma \ref{lem:error-of-x-projected} 
	for the first item in
	\eqref{eq:one-step-projected-weak-error},
	to infer
	\begin{equation}
		\label{eq:auxiliary-projected-euler-weak-error}
		\begin{aligned}
			&
			\quad
			\big| 
			\E \big[
			Z _ E ( t, x; t + h ) -
			Z( t, x; t + h )
			\big]
			\big|
			\\
			&
			=
			\Big|
			\E \Big[
			x - \Phi(x) +    
			\int _ { t } ^ { t + h }
			f(x) - f( \Phi(x) )
			\, \dd s +
			\int _ {t} ^ { t + h }
			g(x) - g(\Phi(x))
			\, \dd W(s)
			\Big]
			\Big|
			%   \\
			%	&
			%	 \leq
			%	  \E [ |
			%	   x - \Phi(x)
			%	    | ]  +    
			%	   \int _ { t } ^ { t + h }
			%	    \E \big[ |
			%	     f(x) - f( \Phi(x) )
			%	      | \big]
			%	     \, \dd s 
			\\
			&
			\leq
			K 
			\big( 1 +
			| x | ^ { 3 \kappa + 4 } 
			\big)
			h ^ { \frac{3}{2} } .
		\end{aligned}
	\end{equation}
	Finally,
	using Lemma \ref{lem:auxiliary-backward-euler-weak-error}
	and
	\eqref{eq:auxiliary-projected-euler-weak-error},
	we have
	\begin{equation}
		\begin{aligned}
			&
			\quad
			\big|
			\E \big[
			X ( t, x; t + h ) -
			Z( t, x; t + h )
			\big]
			\big|
			\\
			&
			=
			\big|
			\E \big[
			( X ( t, x; t + h ) -
			Z _ E ( t, x; t + h ) +
			Z _ E ( t, x; t + h ) -
			Z( t, x; t + h ))
			\big]
			\big|
			\\
			&
			\leq 
			\big|
			\E \big[
			( X ( t, x; t + h ) -
			Z _ E ( t, x; t + h ) )
			\big]
			\big| +
			\big|
			\E \big[
			( Z _ E ( t, x; t + h ) -
			Z( t, x; t + h ) )
			\big]
			\big|
			\\
			&
			\leq
			K 
			\big( 
			1 +
			| x | ^ { 3 \kappa + 4  } 
			\big)
			h ^ { \frac{3}{2} } .
		\end{aligned}
	\end{equation}
	Thus the proof is completed.
	$\square$
	\begin{lem}
		\label{lem:one-step-projected-strong-error}
		Under Assumptions 
		\ref{ass:SDEs-xy-contractive}, \ref{ass:Phi},
		for any 
		$ 1
		\leq
		p 
		\leq
		\frac{ 
  \lfloor
  p ^ {*} 
  \rfloor}
		{ 2\kappa - 1} $,
		we have
		\begin{equation}
			\label{eq:one-step-projected-strong-error}
			\begin{aligned}
				\Big(
				\E \big[
				| X ( t, x; t + h ) -
				Z _E( t, x; t + h ) | ^ {2p}
				\big]
				\Big) ^
				{ \frac{1}{2p} }
				&
				\leq
				K
				\Big
				( 1 +
				| x | ^
				{ ( 4 \kappa - 2 ) p }
				\Big) ^
				{ \frac{1}{2p} } h, 
				\\
				\Big(
				\E \big[
				| Z _ E ( t, x; t + h ) -
				Z( t, x; t + h ) | ^ {2p}
				\big]
				\Big) ^
				{ \frac{1}{2p} }
				&
				\leq
				K
				\Big
				( 1 +
				| x | ^ { 4 (\kappa  + 1 )p + 1 }
				\Big) ^
				{ \frac{1}{2p} } h,
				\\
				\Big(
				\E \big[
				| X ( t, x; t + h ) -
				Z( t, x; t + h ) | ^ {2p}
				\big]
				\Big) ^
				{ \frac{1}{2p} }
				&
				\leq
				K
				\Big
				( 1 +
				| x |^{ 4 (\kappa  + 1 )p + 1 }
				\Big) ^
				{ \frac{1}{2p} } h.
			\end{aligned}
		\end{equation}  
	\end{lem}
	\textit{Proof of Lemma \ref{lem:one-step-projected-strong-error}}.
	First, let us consider the first term in
	\eqref{eq:one-step-projected-strong-error}.
	By using
	\eqref{eq:f-xy-condition},
	\eqref{eq:g-xy-polynomial-condition}, 
	\eqref{eq:X_t-moment-bound}, 
	the moment inequality, 
	the H{\"o}lder inequality and Lemma \ref{lem:one-step-2p-projected-bound},
	one can get 
	\begin{equation}
		\label{eq:auxiliary--euler-strong-error}
		\begin{aligned}
			&
			\quad
			\E \big[
			| X ( t, x; t + h ) -
			Z _E( t, x ; t + h ) | ^ {2p}
			\big]
			\\
			&
			=
			\E \bigg[
			\Big|
			\int _ {t} ^ { t + h }
			f ( X(s) ) -
			f(x)
			\, \dd s +
			\int _ {t} ^ { t + h } 
			g ( X(s) ) -
			g(x)
			\, \dd W(s)
			\Big| ^ {2p}
			\bigg]
			%   \\
			%	& 
			%	 \leq C _ p
			%	  h ^ { 2p - 1 }
			%	   \int _ {t} ^ { t + h }
			%	    \E \big[
			%	     | f ( X(s) ) -
			%	      f(x) | ^ {2p}
			%	       \big]
			%	        \, \dd s +
			%	         C _ p
			%	          h ^ { p - 1 }
			%	          \int _ {t} ^ { t + h }
			%	           \E \big[
			%	            \left\| g ( X(s) ) -
			%	             g (x)
			%	              \right\| ^ {2p}
			%	              \big]
			%	               \,\dd s
			\\
			&
			\leq
			K h ^ { 2p - 1 }
			\int _ {t} ^ { t + h }
			\E \Big[
			\big
			( 1 +
			| X (s) | ^ { 2 \kappa - 2 } +
			|x | ^ { 2 \kappa - 2 }
			\big) ^ {p}
			\big|
			X(s) - x
			\big| ^ {2p}
			\Big]
			\, \dd s 
			\\
			&
			\quad
			+
			K h ^ { p - 1 }
			\int _ {t} ^ { t + h }
			\E \Big[
			\big( 
			1 +
			| X (s) | ^ {  \kappa - 1 } +
			|x | ^ {  \kappa - 1 }
			\big) ^ {p}
			\big|
			X(s) - x
			\big| ^ {2p}
			\Big]
			\, \dd s 
			\\
			& 
			\leq
			K h ^ { p - 1 }
			\int _ {t} ^ { t + h }  
			\Big(
			\E \Big[
			\big
			( 1 +
			| X(s) | ^ { 2 \kappa - 2 } +
			|x | ^ { 2 \kappa - 2 }
			\big) ^ {\frac{(2\kappa-1) p}{\kappa - 1}}
			\Big]
			\Big) ^
			{ \frac{\kappa -1}{2\kappa-1} }
			\\
			&
			\quad            
			\Big(
			\E \Big[
			| X(s) - x | ^ {\frac{2(2\kappa-1)p}{\kappa}}
			\Big]
			\Big) ^
			{ \frac{\kappa}{2\kappa-1} }
			\, \dd s 
			\\
			&
			\leq  
			K 
			\big( 
			1 +
			| x | ^
			{ ( 2 \kappa - 2 ) p } 
			\big)
			( 1 +
			| x | ^
			{ 2 \kappa p } )
			h ^ { p - 1 }
			\int _ {t} ^ { t + h }
			s ^ p
			\, \dd s
			\\
			&
			\leq
			K 
			\big( 1 +
			| x | ^
			{ ( 4 \kappa - 2 ) p } 
			\big)
			h ^ {2p}.
		\end{aligned}
	\end{equation}
	Next,
	utilizing  \eqref{eq:f-xy-condition},
	\eqref{eq:g-xy-polynomial-condition}, 
	\eqref{eq:error-of-x-projected}
	and the  H{\"o}lder inequality,
	the second item of the inequality \eqref{eq:one-step-projected-strong-error} can be estimated as follows:
	\begin{equation}
		\label{eq:auxiliary-projected-strong-error}
		\begin{aligned}
			&
			\quad
			\E \big[
			| ( Z _ E ( t, x; t + h ) -
			Z( t, x; t + h ) ) | ^ {2p}
			\big]
			\\
			&
			=
			\E \bigg[
			\Big|
			x - \Phi(x) +
			\int _ {t} ^ { t + h }
			f (x) -
			f ( \Phi(x) ) 
			\, \dd s +
			\int _ {t} ^ { t + h }
			g(x) -
			g ( \Phi(x) )
			\, \dd W(s)
			\Big| ^ {2p}
			\bigg]
			\\
			&
			\leq
			C _ p
			\E \big[
			| x -
			\Phi(x) | ^ {2p} ] +
			C _ p
			h ^ { 2p - 1 }
			\int _ {t} ^ { t + h }
			\E \big[ 
			| f (x) -
			f ( \Phi(x) ) | ^ { 2p }
			\big]
			\, \dd s 
			\\
			&
			\quad
			+
			C _ p
			h ^ { p - 1 }
			\int _ {t} ^ { t + h }
			\E \big[ 
			| g (x) -
			g ( \Phi(x) ) | ^ { 2p }
			\big]
			\, \dd s
			\\
			&
			\leq
			K 
			\big( 
			1 +
			| x | ^
			{ 4 (\kappa  + 1 )p + 1  } 
			\big)
			h ^ {2p}.         
		\end{aligned}
	\end{equation}
	Further,
	in view of 
	\eqref{eq:auxiliary--euler-strong-error}
	and \eqref{eq:auxiliary-projected-strong-error},
	we obtain
	\begin{equation}
		\begin{aligned}
			&
			\quad
			\E \big[
			| ( X ( t, x; t + h ) -
			Z( t, x; t + h ) ) | ^ {2p}
			\big]
			\\
			&
			=   
			\E \big[
			| ( X ( t, x; t + h ) -
			Z_E ( t, x; t + h ) +
			Z _ E ( t, x; t + h ) -
			Z( t, x; t + h ) ) | ^ {2p}
			\big]
			\\
			& 
			\leq 
			\E \big[
			| ( ( X ( t, x; t + h ) -
			Z_E ( t, x; t + h )) | ^ {2p}
			\big] +
			\E \big[
			| ( Z _ E ( t, x; t + h ) -
			Z( t, x; t + h ) ) | ^ {2p}
			\big]
			\\
			&
			\leq
			K 
			\big( 
			1 +
			| x | ^
			{ 4 (\kappa  + 1 )p + 1  } 
			\big)
			h ^ {2p}.
		\end{aligned}
	\end{equation}
	Thus the proof is finished.
	$\square$

Armed with Theorem \ref{thm:bound-projected-Euler},
Lemmas \ref{lem:one-step-projected-weak-error},
\ref{lem:one-step-projected-strong-error}
and due to Theorem \ref{thm:global-error}, one can straightforwardly obtain the long-time strong convergence rate of the scheme, presented as follows.
	\begin{thm}
		\label{thm:global-error-projected}
		Let Assumptions
		\ref{ass:SDEs-xy-contractive}, \ref{ass:Phi}
  hold.
		For any 
		$ 1
		\leq
		p 
		\leq
		\frac{
  \lfloor
  p ^ {*} 
  \rfloor }
		{ 2 \kappa - 1} $
		and
		$ 0 < h 
		\leq h _ 2$,
		the projected Euler method 
  \eqref{eq:form-modified-Euler}
		has a strong convergence rate of order $\frac{1}{2}$ 
		over infinite time,
		namely,
		\begin{equation}
			\begin{aligned}
				\Big(
				\E \big[
				| X  _ k  -
				Z _ k | ^ {2p}
				\big]
				\Big) ^
				{ \frac{1}{2p} }
				\leq
				K
				\Big(
				\E \big[
				\big( 
				1 +
				| X_0 | ^
				{ ( 6\kappa + 8 ) p }
				\big) 
				\big]
				\Big) ^
				{ \frac{1}{2p} }
				h ^ {\frac{1}{2}}.
			\end{aligned}
		\end{equation}
	\end{thm}

\section{Numerical Experiments}
 \label{sec:numerical experiments}
  \noindent 
In this section, 
we will test the previous findings by performing numerical simulations of some examples of nonlinear SDEs.
\begin{example}
Let us consider the stochastic Ginzburg-Landau (GL) equation \cite{kloeden1992stochastic,hutzenthaler2011strong},
in the form
\begin{equation}
 \label{eq:Ginzburg_Landau_equation}
 \left\{
  \begin{aligned}
   \dd X ( t )
   &
   =
    \big[
     ( \eta +
      \tfrac{1}{2}
       \sigma ^ 2 )
        X ( t ) -
         \vartheta X ^ 3 ( t ) 
          \big]
           \, \dd t
       +
        \sigma 
         X ( t ) 
          \, \dd W ( t ),
           \quad
     t \in [ 0, T ],
           \\
            X(0) 
            &
            =
            X_0>0,
  \end{aligned}
   \right.
\end{equation} 
where 
$ \sigma, 
   \vartheta > 0 $ 
   and $ W : [0,T] \times \Omega \rightarrow \mathbb{R} $ 
   is the real-valued standard Brownian motions.
\end{example}
The coefficients are set as:
$ \eta = - 3/2, 
   \sigma = 1 $
   and
   $ \vartheta = 1 $.
It is easy to see that the coefficients satisfy
Assumptions
\ref{ass:SDEs-xy-contractive} and \ref{ass:g-xy-condition}, 
therefore Theorems \ref{thm:global-error-backward} and
\ref{thm:global-error-projected} are applicable here.
Next,
we will consider the error caused by the temporal discretization of the problem \eqref{eq:Ginzburg_Landau_equation}, 
using both the backward Euler and the projected Euler methods.
In Figure 1 and Figure 2, 
we plot the strong approximation errors of these two numerical schemes for the SDE \eqref{eq:Ginzburg_Landau_equation}.
Set $ T = 16 $ and 
 use the following time step sizes:
$ h 
  \in 
   \{
    2^{-7} , 
     2^{-6} , 
      2^{-5} , 
       2^{-4} , 
        2^{-3} 
          \} 
            $,
where $ h = 2^{-12} $ 
is considered as the exact solution. 
We will use $ M =10000 $ 
 sample paths to simulate the expectation. 
\begin{figure}[!htbp]
	\centering
	 \includegraphics[width=0.6\linewidth, height=0.3
	  \textheight]{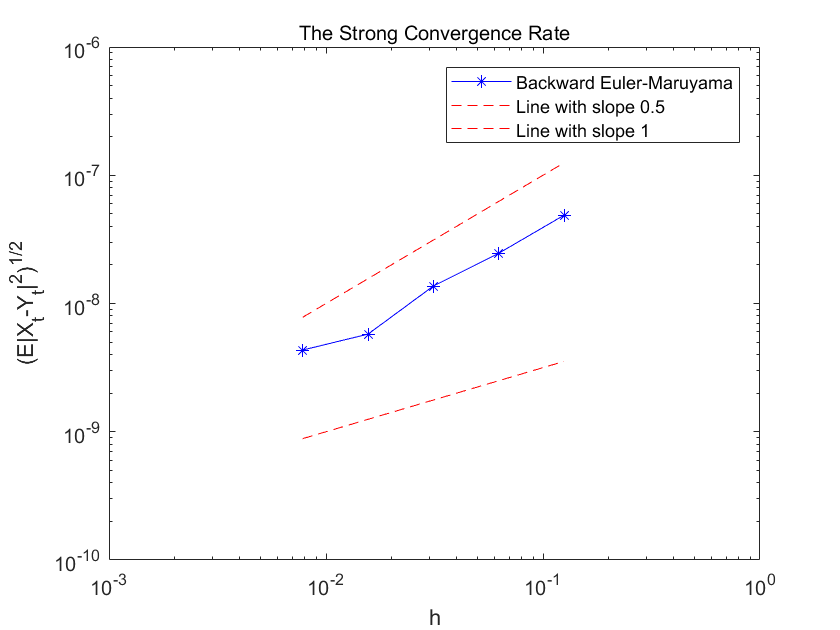}
   	   \caption[Figure 1.]
	{Strong convergence rate of the backward Euler method for
		\eqref{eq:Ginzburg_Landau_equation}.
	}
\end{figure}
\begin{figure}[!htbp]
	\centering
	\includegraphics[width=0.6\linewidth, height=0.3
	\textheight]{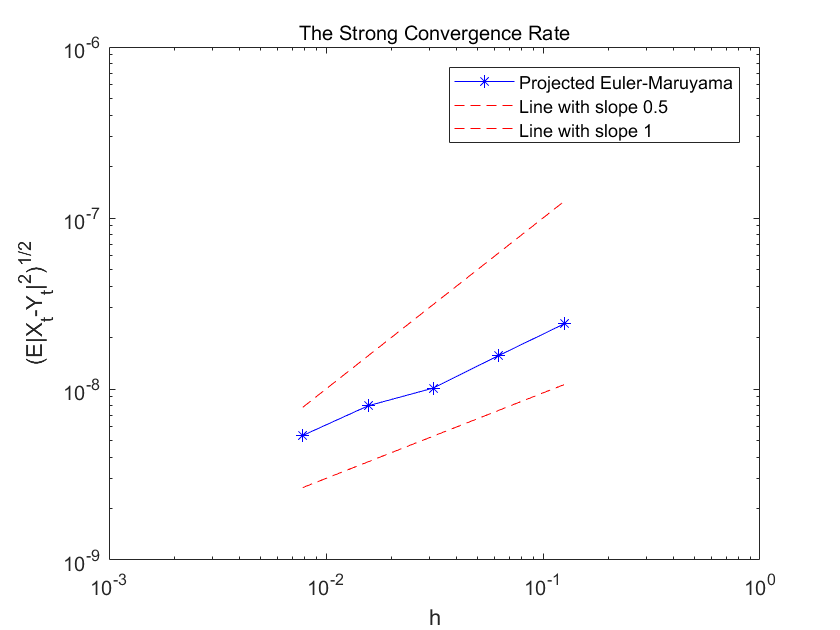}
	\caption[Figure 1.]
	{Strong convergence rate of the projected Euler method for 
		\eqref{eq:Ginzburg_Landau_equation}.
	}
\end{figure}

From Figure 1 and Figure 2, 
the expected strong convergence rate of order 
$ \tfrac{1}{2}$
for both the backward Euler method and
the projected Euler method is numerically confirmed.

\begin{example}
    Consider the following semi-linear stochastic partial differential equation (SPDE):
    \begin{equation}
     \label{eq:SPDE}
	\left\{
          \begin{aligned}      
		\, \dd  
             u(t, x)
            &
             =
             \left[
              \tfrac{\partial^2}{\partial x^2} 
               u(t, x)
        +
         u(t, x)
         -
         u(t, x) ^3
          \right] 
          \, \dd  t
    +
     g(u(t, x)) 
      \, \dd  W_{t }
      ,
      \quad
        t\in (0,T],
        \quad
        x\in(0,1)
      \\
		u(t, 0)
      &
      =
       u(t, 1)
      =
       0 ,
    \\
    u(0, x)
 &
  =
    u_0(x),       
          \end{aligned}
     \right.
\end{equation}
where 
$ g:\mathbb{R}\rightarrow \mathbb{R}$
and
$W:[0,T]\times \Omega \rightarrow \mathbb{R}$ 
is the real-valued standard Brownian motions.
\end{example}
Such an SPDE is usually termed as the stochastic Allen-Cahn
equation.
%To proceed, we aim to spatially discretize the SPDE \eqref{eq:SPDE} in order to derive a system of  SDEs.
We begin by introducing a spatial discretization with a step size 
$\Delta x:=\tfrac{1}{N}$ 
on the interval $[0,1]$ and denoting the discrete spatial points as 
$x_i=i\Delta x,i=1,2,\cdots,N-1$.
 The discretization yields an SDE system:
\begin{equation}
\label{eq:SLSDE5}
 \begin{aligned} 
	\, \dd  X_t
      =
      \left[
       \mathbb{A} 
        X_t
    +
        \mathbb{F}
         (X_t)
          \right]
           \, \dd  t
    +
     \mathbb{G}
      (X_t) 
       \, \dd  W_t, 
        \quad t \in(0, T], 
         \quad X_0=x_0,
 \end{aligned}
\end{equation}
where 
   $X_t=\left(X_{1, t}, X_{2, t}, \cdots, X_{N-1, t}\right)^T:=\left(u\left(t, x_1\right), u\left(t, x_2\right), \cdots, u\left(t, x_{N-1}\right)\right)^T$, 
   $\mathbb{A} \in \mathbb{R}^{(N-1) \times(N-1)}$ ,
   $x_0=\left(u_0\left(x_1\right), u_0\left(x_2\right), \cdots, u_0\left(x_{N-1}\right)\right)^T$ 
and
\begin{center}
    $\mathbb{A}=K^2\left(\begin{array}{cccccc}
	-2 & 1 & 0 & \cdots & 0 & 0 \\
	1 & -2 & 1 & \cdots & 0 & 0 \\
	0 & 1 & -2 & \cdots & 0 & 0 \\
 \vdots	& \vdots &\vdots &   &  \vdots &  \vdots\\
	0 & 0 & 0 & \cdots & -2 & 1 \\
	0 & 0 & 0 & \cdots & 1 & -2
\end{array}\right),$
\end{center}
\begin{center}
    $ \mathbb{F}(X)=\left(\begin{array}{c}
	  f \left(X_1\right)  \\
	 f  \left(X_2\right)  \\ 
	 \vdots \\
	  f\left(X_{N-1}\right) 
     \end{array}\right),\quad \mathbb{G}(X)=\left(\begin{array}{c}
     g\left(X_1\right) \\ 
     g\left(X_2\right) \\
      \vdots \\ 
      g\left(X_{N-1}\right)
  \end{array}\right).
$
\end{center}
Now we turn our attention to the temporal discretization of the SDE system \eqref{eq:SLSDE5},
using the backward Euler and projected Euler methods. 
%Spatial discretization errors are not under consideration for this analysis.
For the following numerical tests, 
we fix $g(u)=\sin{u}+1$, $ T = 30 $ and
and initialize the system with 
$u_0(x)=1$.
%We proceed to conduct our investigations accordingly.

In Figure 3 and Figure 4, 
we plot the strong approximation errors of the two time-stepping schemes for the SDE system \eqref{eq:SLSDE5} with $K=4$.
%To assess the strong convergence rates of these methods, 
We use the following time step-sizes:
$ h 
  \in 
   \{
     \frac{15}{2^{10}} , 
      \frac{15}{2^9}, 
       \frac{15}{2^8}, 
      \frac{15}{2^7} , 
         \frac{15}{2^6}
          \} 
            $
and the numerical approximation with $ h_{\text{exact}} = \frac{15}{2^{12}}  $ 
is identified as the exact solution. 
Moreover, $ M =5000 $  sample paths are used to approximate the expectation.
From Figure 3 and Figure 4, 
one can tell the expected strong convergence rate of order 
$ \tfrac{1}{2}$
for both the backward Euler method and
the projected Euler method, which confirms the theoretical findings.

\begin{figure}[!htbp]
	\centering
	 \includegraphics[width=0.6\linewidth, height=0.3
	  \textheight]{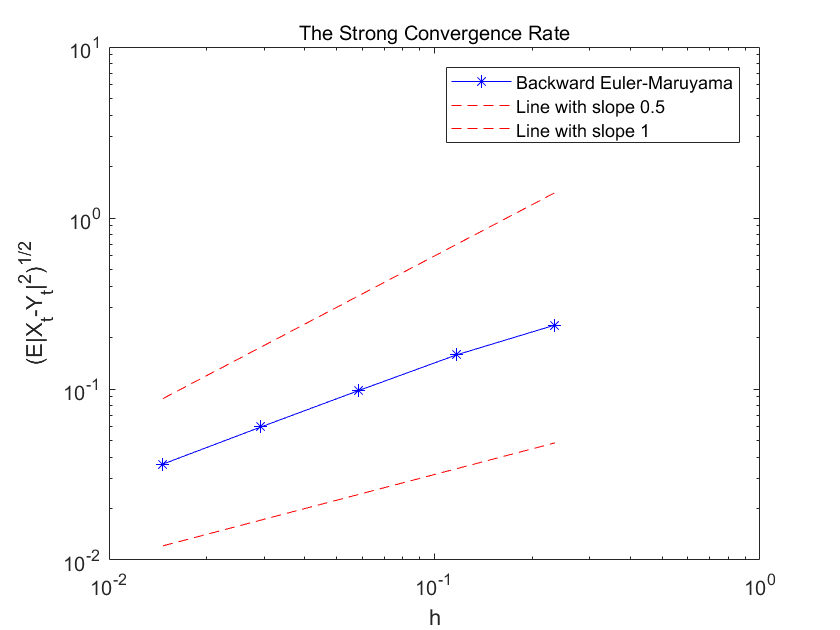}
   	   \caption[Figure 3.]
	{Strong convergence rate of
		the backward Euler method
		for \eqref{eq:SLSDE5}.
	}
\end{figure}
\begin{figure}[!htbp]
	\centering
	 \includegraphics[width=0.6\linewidth, height=0.3
	  \textheight]{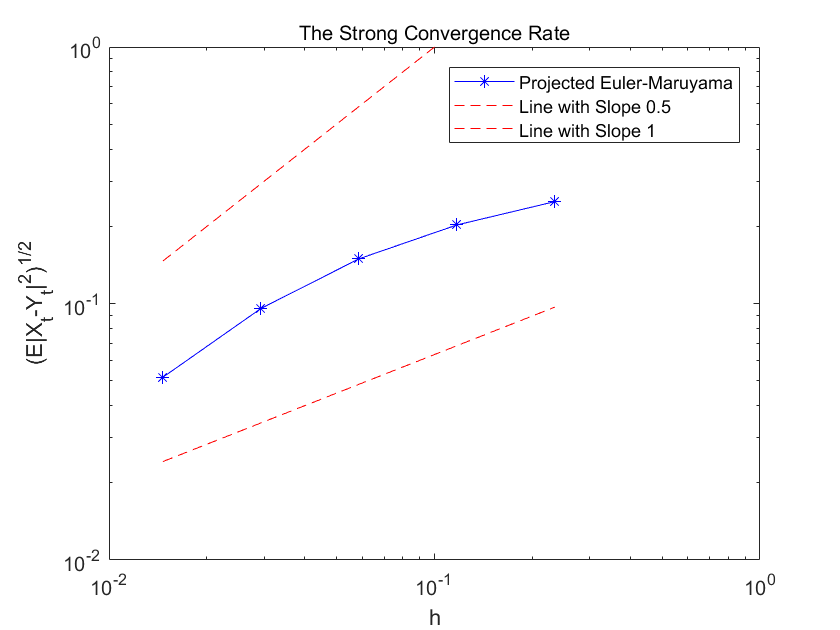}
   	   \caption[Figure 4.]
	{Strong convergence rate of
		the projected Euler method
		for \eqref{eq:SLSDE5}.
	}
\end{figure}

\bibliography{strong_convergence_theorem}

\appendix
\section{Appendix}
\textit{Proof of Theorem 
	\ref{thm:global-error}.}
 To highlight the dependence on initialization, we denote the solution of an SDE as
 $ X ( t _ 0, X _ 0; t _ 0 + t )$.
 Now,
consider the error of the method
  $ Z ( t _ 0, X _ 0; t _ { k + 1 } )$ at the $ ( k + 1 ) $-step
\begin{equation}
 \label{eq:error-decomposition}
  \begin{aligned}
   \rho _ { k + 1 }
  :&
  =
  	X ( t _ 0, X _ 0; t _ { k + 1 } )
  -
  	Z ( t _ 0, X _ 0; t _ { k + 1 } )
  \\
  &
  =
    X ( t _ k, X _ k; t _ { k + 1 } )
  -
    Z ( t _ k, Z _ k; t _ { k + 1 } )
  \\
  &
  =
    X ( t _ k, X _ k; t _ {  k + 1 } ) 
  -
    X ( t _ k, Z _ k; t _ { k + 1 } ) 
  +
    X ( t _ k, Z _ k; t _ { k + 1 } )
  -
    Z ( t _ k, Z _ k; t _ { k + 1 } ).
  \end{aligned}
\end{equation}
It is apparent that the primary distinction on the right-hand side of equation
\eqref{eq:error-decomposition}
arises from the difference on the initial data at time $ t _ k $,
resulting in errors in the solution at the $ ( k + 1 ) $th step.
This can be reformulated as: 
\begin{equation}
   \begin{aligned}
	 \mathcal{D} _ {t _ k, X _ k, Z _ k}( t _ { k+1 } )
    :&
	= 
	  X ( t _ k, X _ k; t _ { k+1 } )
    -
	 X ( t _ k, Z _ k; t _ { k+1 } )
	\\
	&
	:=
	 \rho _ k
	+
	  \mathcal{R}_  {t _ k, X _ k, Z _ k}( t _ { k+1 } )
  	\end{aligned}
\end{equation}
where $  \mathcal{R} $ is given by \eqref{eq:X-xy-form}.  
The second difference in  equality
   \eqref{eq:error-decomposition}
is the one-step error at the
     $ ( k + 1 ) $-step and we denote it as  
\begin{equation}
  	 \mathcal{V} _ { k+1 }
  	:= 
  	X ( t _ k, Z _ k; t _ { k+1 } )
    -
  	Z ( t _ k, Z _ k; t _ { k+1 } ).
\end{equation}
Let $ p \geq 1 $ be an integer.
We have
\begin{equation}
 \label{eq:2p-error-decomposition}
  \begin{aligned}
	\E \big[ |
	 \rho _ { k+1 } |  
	   ^ { 2p } 
	   \big]
  &
  =
    \E \big[
     |  \mathcal{D} _ {t _ k, X _ k, Z _ k}( t _ { k+1 } )
  +
      \mathcal{V} _ { k+1 } | ^ { 2p }
      \big]   
  \\
   &
   =
    \E \big[
     |  \mathcal{D} _ {t _ k, X _ k, Z _ k}( t _ { k+1 } )| ^ { 2 }
   +
     2
     \left\langle
       \mathcal{D} _ {t _ k, X _ k, Z _ k}( t _ { k+1 } ), V _ { k+1 }
       \right\rangle
   + |  \mathcal{V} _ { k+1 } | ^ { 2 }
     \big] ^ { p } 
   \\
	&
    \leq
      \E \big[ 
       |  \mathcal{D} _ {t _ k, X _ k, Z _ k}( t _ { k+1 } ) | ^ { 2p }
        \big]
    +
     2p
      \E \big[
       |  \mathcal{D} _ {t _ k, X _ k, Z _ k}( t _ { k+1 } )| ^ { 2p-2 }
        \left\langle
         \rho _ k +
           \mathcal{R} _ {t _ k, X _ k, Z _ k}( t _ { k+1 } ),  \mathcal{V} _ { k+1 }
           \right\rangle
            \big]
    \\
	 &
	 \quad
	 +
	 \widetilde{ K } _ 1
	  \sum _ { l=2 } ^ { 2p }
	   \E \big[
	    | \mathcal{D} _ {t _ k, X _ k, Z _ k}( t _ { k+1 } )| ^ { 2p-l }
	     |  \mathcal{V} _ { k+1 } | ^ { l }
	      \big],
	\end{aligned}
\end{equation}
where $ \widetilde{ K } _ 1 > 0 $ only depends on $p$.
For the first term on the right-hand side of 
\eqref{eq:2p-error-decomposition},
by \eqref{eq:X-xy-2p} 
we have
\begin{equation}
 \label{eq:S-2p-bound}
  \E \big[
   |  \mathcal{D} _ {t _ k, X _ k, Z _ k}( t _ { k+1 } ) | ^ { 2p }
    \big]
 \leq
   \E \big[
    | \rho _ { k } | ^ { 2p }
     \big]
      \exp( -2p \alpha _ 1 h).
\end{equation}
Next,
we perform a further decomposition of the second term
on the right-hand side of  \eqref{eq:2p-error-decomposition} 
as follows:
\begin{equation}
 \label{eq:SRV-decomposition}
  \begin{aligned}
   &
   \quad
    2p
	\E
     \big[
      | \mathcal{D} _ {t _ k, X _ k, Z _ k}( t _ { k+1 } )| ^ { 2p-2 }
       \left\langle
        \rho _ k +
          \mathcal{R} _ {t _ k, X _ k, Z _ k}( t _ { k+1 } ),  \mathcal{V} _ { k+1 }
         \right\rangle 
          \big]
 \\
  &
  =
   2p
    \E \big[
     | \rho _ { k } | ^ { 2p-2 }
      \left\langle   
       \rho _ { k },  \mathcal{V} _ { k+1 } 
        \right\rangle
         \big]
  +
  2p
   \E \big[
    \big(
     | \mathcal{D} _ {t _ k, X _ k, Z _ k}( t _ { k+1 } ) | ^ { 2p-2 } -
      | \rho _ {k} | ^ { 2p-2} 
      \big)
       \left\langle
        \rho _ {k},  \mathcal{V} _ { k+1 }
         \right\rangle
          \big]
   \\
	&
	\quad 
	+
	2p
	 \E \big[
	  |  \mathcal{D} _ {t _ k, X _ k, Z _ k}( t _ { k+1 } ) | ^ {2p-2}
	   \left\langle
	      \mathcal{R} _ {t _ k, X _ k, Z _ k}( t _ { k+1 } ),  \mathcal{V} _ {k+1}
	      \right\rangle
	       \big].
 \end{aligned}
\end{equation}
Due to the 
$\mathcal{F} _ { t _ k }$-measurability of $\rho _ k$
and the conditional variant of 
\eqref{eq:one-step-weak-error},
we obtain for the first term on the right-hand side of 
\eqref{eq:SRV-decomposition},
\begin{equation}
  2p
  \E \big[
   | \rho _ {k} | ^ { 2p - 2 }
    \left\langle
     \rho _ {k},  \mathcal{V} _ {k+1}
      \right\rangle
       \big]
  \leq
     2p
      C _ 3
      \E \big[
      | \rho _ {k} | ^ { 2p - 1 }
       \big(
        1  +
        | Z _ {k} | ^ { \eta _ 1 } 
        \big) 
         \big] 
          h ^ { q _ 1 }.
\end{equation}
For the second term on the right-hand side of 
\eqref{eq:SRV-decomposition},
noting that it equals zero  when $ p = 1 $,
we have for $ p \geq 2 $,
\begin{equation}
  \begin{aligned}
  &
  \quad	
   2p
   \E
    \big[
     \big( 
      |  \mathcal{D} _ {t _ k, X _ k, Z _ k}( t _ { k+1 } ) | ^ { 2p - 2 }  -
       | \rho _ {k} | ^ { 2p - 2 } 
        \big)
         \left\langle
          \rho _ {k},  \mathcal{V} _ {k+1}
           \right\rangle 
            \big]
   \\
	&
	 \leq
	  2p
	  \widetilde{K} _ 2
	   \E \Big[
	    |  \mathcal{R} _ {t _ k, X _ k, Z _ k}( t _ { k+1 } ) |
	     | \rho _ k |
	      |  \mathcal{V} _ { k + 1 } |
	       \sum _ { l = 0 } ^ { 2p - 3 }
	        |  \mathcal{D} _ {t _ k, X _ k, Z _ k}( t _ { k+1 } ) | ^ { 2p - 3 - l }
	         | \rho _ k | ^ { l } 
	          \Big],
	\end{aligned}
\end{equation}
where $ \widetilde{K} _ 2 > 0 $ only depends on $ p $.
Additionally,
one can utilize $ \mathcal{F} _ { t _ k } $-measurability
of $ \rho _ k $ and the conditional variants of
   \eqref{eq:one-step-strong-error},
    \eqref{eq:X-xy-2p} and \eqref{eq:R-xy-2p},
along with the H{\"o}lder inequality,
to derive for $ p \geq 2 $,
\begin{equation}
 \begin{aligned}
  &
  \quad
   2p
   \E
    \big[
     \big( 
      | \mathcal{D} _ {t _ k, X _ k, Z _ k}( t _ { k+1 } ) | ^ { 2p - 2 } -
       | \rho _ {k} | ^ { 2p - 2 } 
        \big)
         \left\langle
          \rho _ {k},  \mathcal{V} _ { k + 1 }
           \right\rangle
            \big]
 \\
  &
  \leq
   2p
   \widetilde{K} _ 2
    C _ 2
     C _ 4
     \E \Big[
      | \rho _ k | ^ { 2p - 1 }
       \big(
        1 +
        | X _ k | ^ { 2 \kappa - 2 } +
         | Z _ k | ^ { 2 \kappa - 2 }
          \big) ^
           { \frac{1}{4} }
            \big(
             1 +
              | Z _ k | ^ { 2p \eta _ 2 }
               \big) ^
                {\frac{1}{2p} }
                  \Big]
  \\
   &
   \quad
     h ^ { q _ 2  +
      \tfrac{1}{2} } 
       \sum _ { l = 0 } ^ { 2p - 3 }
        \exp
         ( - ( 2p - 3 - l ) \alpha _ 1 h )
  \\
   &
   \leq
    2p
     \widetilde{K} _ 3
      \E \Big[
       | \rho _ k | ^ { 2p - 1 }
        \big(
          1 +
           | X _ k | ^ { 2 \kappa - 2 } +
            | Z _ k | ^ { 2 \kappa - 2 } 
             \big) ^
              { \frac{1}{4} }
               \big(
                1 +
                 | Z _ k | ^ { 2p \eta _ 2 }
                  \big) ^
                  { \frac{1}{2p} }
                    \Big]
                     h ^ { q _ 2 +
                     \frac{1}{2}},
	\end{aligned}
\end{equation}
where $ \widetilde{K} _ 3 > 0 $ depend on $ \widetilde{K} _ 2 $,
$ C _ 2 $ and $C _ 4 $. 

Moreover,
one can bound the third term on the right-hand side of
\eqref{eq:SRV-decomposition},
by employing the conditional variants of
   \eqref{eq:one-step-strong-error}, 
    \eqref{eq:X-xy-2p},
     \eqref{eq:R-xy-2p},
along with applying the H{\"o}lder inequality (twice):
%the term treated in the following way:
  \begin{equation}
   \begin{aligned}
   &
   \quad
    2p
    \E \big[
     |  \mathcal{D} _ {t _ k, X _ k, Z _ k}( t _ { k+1 } )| ^ { 2p - 2 }
      \left\langle
        \mathcal{R} _ {t _ k, X _ k, Z _ k}( t _ { k+1 } ),  \mathcal{V} _ { k + 1 }
       \right\rangle
        \big]
  \\
   &
   \leq
    2p
     \E \Big[
      \Big(
       \E \big[
        |  \mathcal{D} _ {t _ k, X _ k, Z _ k}( t _ { k+1 } )| ^ { 2p }
         \big|
          \mathcal{F} _ { t _ k }
           \big] 
            \Big) ^
             { \frac{2p-2}{2p} }
              \Big(
               \E \big[
                \left\langle
                  \mathcal{R} _ {t _ k, X _ k, Z _ k}( t _ { k+1 } ), V _ { k + 1 } 
                 \right\rangle ^ {p}
                  \big|
                   \mathcal{F} _ { t _ k }
                    \big]
                     \Big) ^
                      { \frac{2}{2p} }
                       \Big]
   \\
    & 
    \leq
     2p
      \E \Big[
       \Big(
        \E \big[
         |  \mathcal{D} _ {t _ k, X _ k, Z _ k}( t _ { k+1 } ) | ^ { 2p }
          \big|
           \mathcal{F} _ { t _ k }
            \big]
             \Big) ^
              { \frac{2p-2}{2p} }
               \Big(
                \E \big[
                 |  \mathcal{R} _ {t _ k, X _ k, Z _ k}( t _ { k+1 } ) | ^ {2p}
                 \big|
                  \mathcal{F} _ { t _ k }
                   \big]
                    \Big) ^
                    { \frac{1}{2p} }
                     \Big(
                      \E \big[
                       |  \mathcal{V} _ { k + 1 } | ^ {2p}
                        \big|
                         \mathcal{F} _ { t _ k }
                          \big]
                           \Big) ^
                            { \frac{1}{2p} }
                             \Big]
  \\
   & 
   \leq
    2p
    C _ 2
     C _ 4
      \E \Big[
       | \rho _ k | ^ { 2p - 1 }
        \big(
          1 +
           | X _ k | ^ { 2 \kappa - 2 } +
            | Z _ k | ^ { 2 \kappa - 2 }
             \big) ^
              { \frac{1}{4} }
               \big(
                1 +
                 | Z _ k | ^ { 2p \eta _ 2 }
                  \big) ^
                   { \frac{1}{2p} }
                    \Big]
                     h ^ { q _ 2 +
                     \frac{1}{2} }
                      \exp
                       \big(
                        - ( 2p - 2 ) \alpha _ 1 h
                         \big) 
  \\
   &
   \leq
    2p
     \widetilde{K} _ 4
      \E \Big[
       | \rho _ k | ^ { 2p - 1 }
        \big(
         1 +
           | X _ k | ^ { 2 \kappa - 2 } +
             | Z _ k | ^ { 2 \kappa - 2 }
              \big) ^
              { \frac{1}{4} }
               \big(
                 1 +
                 | Z _ k | ^ { 2p \eta _ 2 }
                  \big) ^
                   { \frac{1}{2p} }
                    \Big]
                     h ^ { q _ 2 +
                      \frac{1}{2} },
 \end{aligned}
\end{equation}
where $ \widetilde{K} _ 4 > 0 $ depends on $ C _ 2, C _ 4 $.  

Following the conditional versions of
    \eqref{eq:one-step-strong-error}
and \eqref{eq:X-xy-2p}
as well as the H{\"o}lder inequality,
one can bound the third  third term on the right-hand side of the inequality \eqref{eq:2p-error-decomposition}  as follows:
  \begin{equation}
   \label{eq:SV-error}
  	\begin{aligned}
  &
  	\quad
  	 \widetilde{K} _ 1
  	  \sum _ { l = 2 } ^ { 2p }
  	   \E \big[
  	    |  \mathcal{D} _ {t _ k, X _ k, Z _ k}( t _ { k+1 } ) | ^ { 2p - l }|
  	      \mathcal{V} _ { k + 1 }| ^ {l}
  	      \big]
  \\
   & 
   \leq  
    \widetilde{K} _ 1
     \sum _ { l = 2 } ^ { 2p }
      \E \Big[
       \Big(
        \E \big[
        |  \mathcal{D} _ {t _ k, X _ k, Z _ k}( t _ { k+1 } ) | ^ { 2p }
         \big|
          \mathcal{F} _ { t _ k }
           \big]
            \Big) ^
             { \frac{ 2p - l}{2p} }
              \Big(
               \E \big[
                |  \mathcal{V} _ { k + 1 } | ^ { 2p }
                \big|
                 \mathcal{F} _ { t _ k }
                  \big]
                   \Big) ^
                    { \frac{l}{2p} }
                     \Big]
  \\
   &
   \leq 
    \widetilde{K} _ 1
     \sum _ { l = 2 } ^ { 2p }
     ( C _ { 4 }) ^ { l }
      \E \Big[
       | \rho _ k | ^ { 2p - l }
        \big(
         1 +
         | Z _ k | ^ { 2p \eta _ 2 }
          \big) ^
           { \frac{l}{2p} }
            \Big]
             h^{ l q _ 2 }
             \exp( - ( 2p - l ) \alpha _ 1 h )
  \\
   &
   \leq
    \widetilde{K} _ 5 
     \sum _ { l = 2 } ^ {2p}
      \E \Big[
       | \rho _ k | ^ { 2p - l }
        \big(
         1 +
          | Z _ k | ^ { 2p \eta _ 2 }
           \big) ^
            { \frac{l}{2p} }
             \Big]
               h^{ l q _ 2 }, 
  \end{aligned}
\end{equation}
  where $ \widetilde{K} _ 5 > 0 $ depends on
   $ \widetilde{K} _ 1, C _ 4 $.  
By Substituting the inequalities
 \eqref{eq:S-2p-bound} to
 \eqref{eq:SV-error}
into 
   \eqref{eq:2p-error-decomposition}
and recalling that
     $ q _ 1 \geq q _ 2 + \tfrac{1}{2} $,
one can obtain
\begin{equation}
 \label{eq:one-step-2p-error}
  \begin{aligned}
   \E \big[ | 
    \rho _ { k + 1 } 
     | ^ {2p} 
      \big]
  &
  \leq
   \E \big[
    | \rho _ k | ^ {2p}
     \big]
      \exp( -2p \alpha _ 1 h )
  +
   2p
    C _ 3
    \E \big[
     | \rho _ {k}| ^ { 2p - 1 }
      \big( 1 +
       | Z _ {k}| ^ { \eta _ 1 } 
        \big)
        \big]
         h ^ { q _ 1 }
 \\
  &
  \quad
   +
   2p
   \widetilde{K} _ 3
    \E \Big[
     | \rho _ k | ^ { 2p - 1 }
      \big(
       1 +
        | X _ k | ^ { 2 \kappa - 2 } +
         | Z _ k | ^ { 2 \kappa - 2 }
          \big) ^
           { \frac{1}{4} }
             h ^ { q _ 2 +
              \tfrac{1}{2} }
               \big(
                 1 +
                  | Z _ k | ^ { 2p \eta _ 2 }
                   \big) ^
                    { \frac{1}{2p} }
                     \Big]
 \\
  &
  \quad
   +
   2p
   \widetilde{K} _ 4
    \E \Big[
     | \rho _ k | ^ { 2p - 1 }
      \big(
        1 +
         | X _ k | ^ { 2 \kappa - 2 } +
          | Z _ k | ^ { 2 \kappa - 2 }
           \big) ^
            { \frac{1}{4} }
              h ^ { q _ 2 +
              \frac{1}{2} }
               \big(
                1 +
                | Z _ k | ^ { 2p \eta _ 2 }
                 \big) ^
                  { \frac{1}{2p} }
                   \Big] 
 \\
  &
  \quad
   +
   \widetilde{K} _ 5
    \sum _ { l = 2 } ^ {2p}
     \E \Big[
       | \rho _ k | ^ { 2p - l }
        h ^ { l q _2 }
        \big(
         1 +
          | Z _ k | ^ { 2p \eta _ 2 }
           \big) ^
            { \frac{l}{2p} }
             \Big]. 
% \\
%  & 
%  \leq
%   \E \big[
%   | \rho _ k | ^ {2p}
%     \big]
%      \exp( - 2p \alpha _ 1 h )
%  +
%  \frac{ 2p \alpha _ 1 h}{16}
%   \E \Big[
%    | \rho _ {k}  | ^ { 2p - 1 }
%     \widetilde{K}
%      \big(
%       1 +
%       | Z _ {k}| ^ { \eta _ 1 }
%        \big)
%         h ^ { q _ 2 -
%         \tfrac{1}{2} }
%          \Big] 
% \\
%  &
%  \quad
%   +
%   \frac{ 2p \alpha _ 1 h }{16}
%    \E \Big[
%     | \rho _ {k}  | ^ { 2p - 1 }
%      \widetilde{K}
%       \big(
%        1 +
%         | X _ k | ^ { 2 \kappa - 2 } +
%          | Z _ k| ^ { 2 \kappa - 2 }
%           \big) ^
%            { \tfrac{1}{4} }
%              h ^ { q _ 2 -
%               \tfrac{1}{2} }
%                \big(
%                  1 +
%                   | Z _ k | ^ { 2p \eta _ 2 }
%                    \big) ^
%                     { \tfrac{1}{2p} }
%                       \Big]
% \\
%  &
%  \quad 
%   +
%   \frac{ \alpha _ 1 h }{8}
%    \sum _ { l = 2 } ^ {2p}
%     \E \Big[
%      | \rho _ k | ^ { 2p - l }
%       \widetilde{K}
%        \big(
%         1 +
%          | Z _ k | ^ { 2p \eta _ 2 }
%           \big) ^
%            { \tfrac{l}{2p} }
%             h^{ l q _ 2 - 1 }
 %             \Big],     
    \end{aligned}
     \end{equation}
Then using  the Young inequality for 
 \eqref{eq:one-step-2p-error},
we have
\begin{equation}
 \begin{aligned}
  \E \big[
   | \rho _ { k + 1 } | ^ {2p}
    \big]
  &
   \leq
    \E \big[ | 
     \rho _ k
       | ^ {2p} 
        \big]
         \exp( - 2p \alpha _ 1 h )
  +
   \frac{ p \alpha _ 1 h }{8} 
    \E \big[ |
     \rho _ {k}
      | ^ {2p} 
       \big]
  +
  \widetilde{K} 
   \E \big[
    \big(
     1 +
      | Z _ k | ^ { \eta _ 1 }
       \big) ^ {2p}
        \big]
         h ^ { 2p
      	  ( q _ 2 -
          \frac{1}{2} ) +
           1  } 
 \\
  &
  \quad
   +
   \frac{ p \alpha _ 1 h }{8} 
     \E \big[ | 
      \rho _ {k}
       | ^ {2p}
        \big]
  + 
  \widetilde{K}
   \E \Big[
    \big(
     1 +
      | X _ k | ^ { 2 \kappa - 2 } +
       | Z _ k | ^ { 2 \kappa - 2 }
        \big) ^
         { \frac{p}{2} }
          \big( 1 +
           | Z _ k | ^ { 2p \eta _ 2 }
            \big) 
             \Big]
              h ^ { 2p ( q _ 2 -
              \frac{1}{2} ) +
                1 } 
 \\
  &
  \quad 
   +
   \frac{
   	\sum _ { l = 2 } ^ {2p}
   	 \tfrac{2p-l}{2p} } {8}
      \alpha _ 1 h
       \E \big[ | 
        \rho _ k
         | ^ {2p} 
          \big]
   +
   \widetilde{K}
    \E \big[
     \big(
      1  +
       | Z _ k | ^ { 2p \eta _ 2 }
        \big) 
         \big]
          h ^ { 2p
       	   ( q _ 2 -
            \frac{1}{2} ) +
              1} ,  	
   \end{aligned}
\end{equation}
where $  \widetilde{K} > 0 $ is independent of $ h, t $.

Given the assumption that 
 $ 0 < h \leq \frac{1}{p \alpha _ 1} $
 and the inequality 
 $ e ^ {-x} \leq 1 - x + \frac { x ^ {2} }{2} $
for $ 0 < x < 1 $,
we can proceed by utilizing the inequalities
     $\sum _ { l = 2 } ^ {2p} \frac { 2p - l }{2p} \leq p $, 
     \eqref{eq:X_t-moment-bound} and 
     \eqref{eq:numerical-solution-bound} and obtain
\begin{equation}
 \begin{aligned}
  \E \big[
   | \rho _ { k + 1 } | ^ {2p}
    \big]
  &
   \leq
    \big( 1 -
     \frac{p}{2}
      \alpha _ 1 h +
       \frac{2p}{8} 
        \alpha _ 1 h +
         \frac{p}{8}
          \alpha _ 1 h 
           \big)
            \E \big[ | 
             \rho _ k
              | ^ {2p} 
               \big] +
                \widetilde{K} 
                 \E \big[
                  \big(
                   1 +
                   | Z _ k | ^ {2p \eta _ 1 }
                    \big) 
                     \big]
                      h ^ { 2p
               	      ( q _ 2 -
               	       \frac{1}{2} ) +
               	         1  } 
 \\
  &
  \quad
   +
   \widetilde{K} 
    \E \Big[
     \big(
      1 +
       | X _ k | ^ { 2 \kappa - 2 } +
        | Z _ k | ^ { 2 \kappa - 2 }
         \big) ^
         { \tfrac{p}{2} }
          \big(
            1 +
             | Z _ k | ^ { 2 p\eta _ 2 }
             \big) 
              \Big]
               h^{ 2p
               	( q _ 2 -
                 \frac{1}{2} ) +
                  1}
 \\
  &
   \leq
    \big( 1 -
     \frac{p}{8} 
      \alpha _ 1 h 
       \big)
        \E \big[
         | \rho _ k | ^ {2p}
          \big]
  +
   \widetilde{K}
    \E \big[ 
     \big(
      1 +
       | X _ 0 | ^ 
        { 2 \lambda p}
         \big)
          \big]
           h ^ { 2p
            ( q _ 2 -
             \frac{1}{2} ) +
              1},
	\end{aligned}
\end{equation}
where 
$ 
\lambda 
= 
\max 
\{ 
\eta_1 \eta_3, 
\big( 
\tfrac{ \kappa - 1 }{2} +
\eta_2 
\big)
\eta_3
 \}. $ 
Equipped with the above inequality,
one sees
\begin{equation}
	\E \big[ | 
	 \rho _ k
	  | ^ {2p} 
	   \big] 
   \leq 
	 C
	  \E \big[ 
	   \big(
	    1 +
	     | X _ 0 | ^ { 2 \lambda p}
	      \big)
	       \big]
	        h ^ { 2p 
	         (q _ 2 -
	          \frac{1}{2} ) },
\end{equation}
where the constant $ C $ does not depend on $h,t,k$.
Namely,
\begin{equation}
 \begin{aligned}
  \big(
   \E \big[
    | X ( t _ 0, X _ 0; t _ k ) -
     Z ( t _ 0, X _ 0; t _ k ) | ^ {2p}
     \big]
      \big) ^
       { \frac{1}{2p} }
  \leq 
     C
     \big(
      \E \big[
       \big(
        1 +
         | X _ 0 | ^ { 2 \lambda p }
          \big) 
            \big]
             \big) ^
             { \frac{1}{2p} }
              h ^ { q _ 2 -
               \frac{1}{2}  }.
	\end{aligned}
\end{equation}
To conclude,
Theorem \ref{thm:global-error} is proved for integer $p$.
For non-integer $p$, the conclusion  can be derived using the Young inequality. Therefore, 
 the proof is completed.
$\square$

\end{document}